\documentclass[reqno, 12pt]{amsart}
\usepackage{dsfont}
\usepackage{amssymb}
\usepackage{mathrsfs}
\usepackage{enumitem}
\usepackage{comment}
\usepackage[backend=biber, style=alphabetic, url=false]{biblatex}
\usepackage[pdfstartview=FitH, pdfborder={0 0 0}, colorlinks=true, citecolor=blue, linkcolor=blue]{hyperref}
\usepackage{aliascnt}

\theoremstyle{plain}
\newtheorem{theorem}{Theorem}[section]

\newaliascnt{conjecture}{theorem}

\aliascntresetthe{conjecture}

\newaliascnt{corollary}{theorem}

\aliascntresetthe{corollary}

\newaliascnt{lemma}{theorem}
\newtheorem{lemma}[lemma]{Lemma}
\aliascntresetthe{lemma}

\newaliascnt{fact}{theorem}

\aliascntresetthe{fact}

\newaliascnt{claim}{theorem}

\aliascntresetthe{claim}

\newaliascnt{proposition}{theorem}
\newtheorem{proposition}[proposition]{Proposition}
\aliascntresetthe{proposition}

\theoremstyle{definition}

\newaliascnt{definition}{theorem}

\aliascntresetthe{definition}

\newaliascnt{example}{theorem}

\aliascntresetthe{example}

\theoremstyle{remark}

\newaliascnt{remark}{theorem}
\newtheorem{remark}[remark]{Remark}
\aliascntresetthe{remark}

\numberwithin{equation}{section}

\usepackage{graphicx}
\usepackage{tikz}
\usetikzlibrary{decorations.pathreplacing}
\allowdisplaybreaks[3]
\theoremstyle{plain}
\newtheorem{theoremx}{Theorem}

\newcommand{\Sph}{\mathbb S}
\newcommand{\Ric}{\operatorname{Ric}}
\newcommand{\rank}{\operatorname{rank}}

\begin{document}

\title[Hemisphere Rigidity]{Saddle Codazzi tensor on spheres and hemisphere rigidity}

\author{Jian Ge}
\address{School of Mathematical Sciences, Laboratory of Mathematics and Complex Systems, Beijing Normal University, Beijing 100875, P. R. China}
\email{jge@bnu.edu.cn}
\thanks{The first-named author is supported in part by NSFC No.12371049 and the Fundamental Research Funds for the Central Universities.}

\author{Mijia Lai}
\address{School of Mathematical Sciences, Shanghai Jiao Tong University,\newline
Shanghai 200240, China}
\email{laimijia@sjtu.edu.cn}

\author{Yuntian Zhu}
\address{School of Mathematical Sciences, Shanghai Jiao Tong University,\newline
Shanghai 200240, China}
\email{fxy012805@sjtu.edu.cn}

\date{}

\begin{abstract}
Let $(M^n,g)$, $n\geq 3$, be a compact connected Riemannian manifold
with smooth boundary. We ask whether $M$ must be isometric to the
standard hemisphere if $\sec_g\geq 1$, $\partial M $ is isometric to
the standard unit sphere, and its second fundamental form is Codazzi.

Our approach centers on saddle Codazzi tensors on round spheres,
where the saddle condition requires the restriction to every tangent
two-plane to have nonpositive determinant.
By analyzing the associated support maps, we prove that every smooth
saddle Codazzi tensor on $\mathbb S^m$ vanishes identically for
$m\geq 3$. Together with Panina's examples on $\mathbb S^2$, this
shows that, among unit spheres of dimension at least two, nontrivial
smooth saddle Codazzi tensors exist precisely in dimension two.
The vanishing theorem yields an affirmative answer to the hemisphere
rigidity question for $n\geq 4$.

In dimension three, we realize the nontrivial saddle Codazzi tensors
provided by Panina's construction, after scaling, as the second
fundamental forms of smooth perturbations of the standard hemisphere,
preserving the round boundary metric and the sectional curvature
lower bound. This produces counterexamples to rigidity in dimension
three. Nevertheless, under the same curvature and boundary assumptions,
we prove the sharp inradius estimate $r_{\mathrm{in}}(M)\leq \pi/2$,
with equality if and only if $M$ is isometric to the standard hemisphere.
\end{abstract}
\maketitle

\section{Introduction}

Rigidity phenomena for compact Riemannian manifolds with boundary have been studied extensively from both geometric and analytic points of view. A fundamental result in this direction is a theorem by Hang and Wang, \cite{HW}. In one of its standard forms, their theorems asserts that if $(M^n, g)$ is a compact Riemannian manifold with smooth boundary such that
\begin{equation*}
	\Ric_g\geq (n-1)g, \qquad (\partial M,g|_{T\partial M})\cong\Sph^{n-1}, \qquad h\geq0,
\end{equation*}
where $h$ denotes the second fundamental form of $\partial M$ with respect to the outward unit normal and the convention $h(X,Y)=g(\nabla^M_X\nu,Y)$, then $(M, g)$ is isometric to the standard hemisphere $\Sph^n_+$. Generalizations to spherical domains bounded by Clifford tori appear in~\cite{L}.

The prescribed boundary metric and second fundamental form can also be viewed as geometric analogues of Dirichlet and Neumann
data. This viewpoint connects hemisphere rigidity with overdetermined boundary problems, such as Serrin's symmetry
theorem~\cite{S}, and with rigidity results for the Obata equation $\nabla^2f+fg=0$; see~\cite{R,E,X,CLW}.
We ask whether the stronger interior condition $\sec_g\geq1$ allows one to replace boundary convexity by a condition that
does not prescribe the signs of the principal curvatures.

Grove and Petersen~\cite{GP} studied the related problem for
Alexandrov spaces of curvature at least $1$. They prove that
maximal boundary volume,
$\operatorname{vol}(\partial M)=\operatorname{vol}(\Sph^{n-1})$,
forces the space to be a spherical lens $L^n_\alpha$, bounded by
two hemispheres meeting at an angle $\alpha\in(0,\pi]$.
This does not directly resolve the smooth Riemannian problem:
without local convexity of the boundary, a lower sectional
curvature bound need not make the intrinsic distance space an
Alexandrov space with the same lower curvature bound.

Our additional boundary assumption is the Codazzi equation.
After identifying the boundary with the unit sphere, it reads
\begin{equation*}
	(\nabla_Xh)(Y,Z)=(\nabla_Yh)(X,Z)
	\qquad (X,Y,Z\in T\Sph^{n-1}).
\end{equation*}
Together with the Gauss equation, this condition gives a
boundary problem on the round sphere that can be studied
independently of the interior metric. The resulting rigidity
holds in dimensions at least four.

\begin{theoremx}\label{thm:main}
	Let $(M^n,g)$, $n\geq4$, be a compact connected Riemannian
	manifold with smooth boundary. Assume that
	\begin{enumerate}
		\item $\sec_g\geq1$;
		\item $(\partial M,g|_{T\partial M})$ is isometric to the
		standard unit sphere $\Sph^{n-1}$;
		\item the second fundamental form $h$ of $\partial M$ is a
		Codazzi tensor.
	\end{enumerate}
	Then $(M,g)$ is isometric to the standard hemisphere $\Sph^n_+$.
\end{theoremx}

\begin{remark}
	The boundary argument uses only $\sec_M(\Pi)\geq1$ for planes
	$\Pi\subset T\partial M$. Once $h=0$, the theorem of Hang and
	Wang uses $\Ric_g\geq(n-1)g$ in the interior. These two
	assumptions may therefore replace condition~(1) in
	Theorem~\ref{thm:main}.
\end{remark}

\begin{remark}
	The corresponding statement is false in dimension two. A
	smooth simple closed curve of length $2\pi$ in $\Sph^2$ need
	not be a great circle. Either region bounded by such a curve
	has Gaussian curvature $1$ and boundary intrinsically
	isometric to $\Sph^1$, but is not a hemisphere. The Codazzi
	condition imposes no restriction on a one-dimensional boundary.
\end{remark}

Dimension three admits a different behavior. Panina's
construction~\cite{Pan05} provides nonzero Codazzi tensors on
$\Sph^2$ satisfying the sign condition required by the Gauss
equation. We show that these boundary data can be realized by
metrics whose sectional curvatures remain at least $1$
throughout the interior.

\begin{theoremx}\label{thm:three}
	There exists a compact Riemannian manifold $(M^3,g)$,
	diffeomorphic to the closed three-dimensional ball, such that
	\begin{equation*}
		\sec_g\geq1,
		\qquad
		(\partial M,g|_{T\partial M})\cong\Sph^2,
	\end{equation*}
	and the second fundamental form of $\partial M$ is a nonzero
	Codazzi tensor. In particular, $(M,g)$ is not isometric to the
	standard hemisphere.
\end{theoremx}

Theorem~B is related to the study of saddle spheres in the standard
$\Sph^3$. G\'alvez, Mira, and Tassi \cite{GMT24} prove that every
real analytic immersed saddle two-sphere in $\Sph^3$ is equatorial,
and show that this rigidity fails in the smooth category.
In their setting, the ambient metric is round, while the induced
metric on the sphere is allowed to vary. Our construction fixes
the induced boundary metric to be the standard metric on $\Sph^2$
and realizes nontrivial saddle Codazzi tensors, after scaling,
as its second fundamental form by deforming the ambient metric
while maintaining $\sec_g\geq1$.

Although these boundary conditions do not imply rigidity in
dimension three, they give a sharp bound for the distance from
an interior point to the boundary. Equality restores rigidity.

\begin{theoremx}\label{thm:inradius-three}
	Let $(M^3,g)$ be a compact connected Riemannian manifold with
	smooth boundary. Assume that
	\begin{enumerate}
		\item $\sec_g\geq1$;
		\item $(\partial M,g|_{T\partial M})$ is isometric to the
		standard unit sphere $\Sph^2$;
		\item the second fundamental form $h$ of $\partial M$ is a
		Codazzi tensor.
	\end{enumerate}
	Then
	\begin{equation*}
		r_{\mathrm{in}}(M)
		:=\max_{p\in M}d_g(p,\partial M)\leq\frac{\pi}{2},
	\end{equation*}
	with equality if and only if $(M,g)$ is isometric to the
	standard hemisphere $\Sph^3_+$.
\end{theoremx}

The common geometric mechanism behind our proofs is the study of
\emph{saddle Codazzi tensors} on round spheres. These are symmetric
Codazzi tensors whose restrictions to all tangent two-planes have
nonpositive determinant. Put $m=n-1$ and identify a symmetric
boundary tensor with its self-adjoint endomorphism using the round
metric. Every Codazzi tensor on $\Sph^m$, $m\geq2$, is of the form
$A_f=\nabla^2f+fg_{\Sph^m}$. When the boundary is round and
$\sec_g\geq1$, the Gauss equation shows that its second fundamental
form satisfies precisely the saddle condition. The boundary geometry
thus leads naturally to the question of whether nontrivial saddle
Codazzi tensors can exist on $\Sph^m$. Pointwise, the saddle condition
allows at most one positive and one negative eigenvalue, so
$\rank A_f\leq2$. This rank restriction is the starting point of
our analysis of the associated support map.

The potential $f$ determines a map into Euclidean space,
\begin{equation*}
	F_f:\Sph^m\longrightarrow\mathbb R^{m+1},
	\qquad F_f(x)=\nabla f(x)+f(x)x,
	\qquad dF_f=A_f.
\end{equation*}
It parametrizes the envelope of the planes
$\langle y,x\rangle=f(x)$, with singularities allowed.
On a sufficiently small constant-rank neighborhood its image
is a submanifold, and $x$ is a unit normal to that image.
If $r=\max\rank A_f>0$ and $m\geq3$, then $r\leq2<m$.
A suitable linear projection of $F_f$ has a compact regular
fiber component that is a full unit normal sphere of a local
image. This sphere has dimension $m-r\geq1$, so it contains
a complete great circle on which $F_f$ is constant.

Such a circle gives a global obstruction. Restrict $f$ to any
great $2$-sphere containing the circle and subtract the linear
function determined by the constant value of $F_f$. The resulting
function $u$ and its first derivatives vanish along the circle.
Vanishing of $u$ alone makes the orthogonal projection of $F_u$
onto the plane of the circle zero there. Stokes' theorem gives
zero projected signed area on each hemisphere, so the saddle
condition forces $\det A_u=0$. An energy identity then realizes
equality in the first Dirichlet eigenvalue inequality on each
hemisphere and shows that $u$ is linear. Its vanishing first
derivatives along the circle now give $u=0$.
Varying the great $2$-sphere gives $A_f=0$ on $\Sph^m$.
Thus $h=0$ in every
dimension $n\geq4$, and Hang--Wang rigidity applies.

The same map also explains the boundary data in dimension
three. Its signed area vector gives the balance identity needed
to construct the counterexamples. We first realize the data on
an exact boundary collar. A correction supported in the interior
then gives the global curvature bound while preserving the collar.

For the inradius equality case, we separate distance comparison
from a rigidity statement for contact sets on $\Sph^2$
(Proposition~\ref{prop:three-contact-rigidity}). Distance comparison
produces a nonnegative boundary function $u$ satisfying
$\nabla^2u+b h\leq0$, with $b>0$, near its compact zero set $C$.
At contact points, $h\leq0$ and $\rank h\leq1$.
The two distinct support-map values supplied by
G\'alvez--Mira~\cite{GM}, together with a convex barrier argument,
force a scalar projection of $F_f$ to have an interval of values
on $C$. A regular level curve meeting $C$ stays in $C$ and becomes
a complete great circle on which $F_f$ is constant.
Lemma~\ref{lem:circle-support-rigidity} then gives $h=0$.
Thus the two rigidity arguments end with the same circle
obstruction; in dimension three, the contact inequality and the
G\'alvez--Mira theorem supply the circle.

Section~\ref{sec:support-geometry} develops support maps and
the geometry of their normal fibers, which yields non-existence of nontrivial saddle Codazzi tensors on $\Sph^m$ for $m\geq 3$. Using this, we complete the proof of
Theorem~\ref{thm:main}.
Section~\ref{sec:counterexamples} realizes the nonzero saddle
data on $\Sph^2$ by metrics with $\sec\geq1$.
Section~\ref{sec:dimension-three} proves the inradius estimate
and its equality case.

\section{Classification of Saddle Codazzi tensors on $\Sph^m (m\geq3)$ }\label{sec:support-geometry}

On the unit round sphere $(\Sph^m,g_{\Sph^m})$, $m\geq2$, we identify
a symmetric bilinear form $h$ with its self-adjoint endomorphism $A$
by $h(X,Y)=\langle AX,Y\rangle$. We say that $A$ is \emph{saddle} if
\[
	\det(h|_\Pi)\leq0
	\qquad\text{for every two-plane }\Pi\subset T_x\Sph^m
	\text{ and every }x\in\Sph^m,
\]
where the determinant is computed in an orthonormal basis of $\Pi$.
For a function $f$, we write $A_f$ for the endomorphism determined by
\begin{equation*}
	\langle A_fX,Y\rangle
	=\nabla^2f(X,Y)+f\langle X,Y\rangle.
\end{equation*}
Thus $A_f$ represents the spherical Hessian of $f$ under this
identification. Pointwise evaluation is denoted by $A(x)$ or $A_f(x)$;
the subscript in $A_f$ records the potential function.
For the boundary in Theorem~\ref{thm:main}, the Gauss equation gives
\begin{equation}\label{eq:boundary-saddle-gauss}
	\det(h|_\Pi)=1-\sec_M(\Pi)\leq0.
\end{equation}
The essential rigidity statement is therefore intrinsic to the sphere.

We will prove that this boundary tensor vanishes in every dimension
$m\geq3$. We first describe the support map and the geometry of its compact fibers.

\begin{lemma}\label{lem:rank-two}
	A self-adjoint endomorphism is saddle if and only if it has at most
	one positive and at most one negative eigenvalue. In particular,
	its rank is at most two.
\end{lemma}

\begin{proof}
	Applying the saddle condition to the plane spanned by any two
	orthonormal eigenvectors gives $\lambda_i\lambda_j\leq0$ for $i\neq j$.
	Conversely, if there is at most one eigenvalue of each sign, the
	restriction to any two-plane is neither positive nor negative
	definite, and hence has nonpositive determinant.
\end{proof}

\begin{lemma}\label{lem:codazzi-potential-map}
	Let $A$ be a smooth self-adjoint Codazzi tensor on $\Sph^m$,
	$m\geq2$. Then $A=A_f$ for a smooth function $f$, and the map
	\begin{equation*}
		F:\Sph^m\longrightarrow\mathbb R^{m+1},
		\qquad F(x)=\nabla f(x)+f(x)x,
	\end{equation*}
	satisfies $dF_x(X)=A(x)X$, viewed as an identity in
	$\mathbb R^{m+1}$. Conversely, $A_f$ is Codazzi for every smooth $f$.
	The potential $f$ is unique modulo functions
	$x\mapsto\langle a,x\rangle$.
\end{lemma}

\begin{proof}
	Regard $\omega(X)=AX$ as an $\mathbb R^{m+1}$-valued one-form.
	The Euclidean Gauss formula gives
	\begin{equation*}
		d\omega(X,Y)
		=(\nabla_XA)Y-(\nabla_YA)X
		-\bigl(\langle X,AY\rangle-\langle Y,AX\rangle\bigr)x=0.
	\end{equation*}
	Codazzi cancels the tangential part, and self-adjointness cancels
	the normal part. Since $\Sph^m$ is simply connected, $\omega=dF$
	for a smooth map $F:\Sph^m\to\mathbb R^{m+1}$. Put
	$f(x)=\langle F(x),x\rangle$. As $dF(X)=AX$ is tangent to the
	sphere, $df(X)=\langle F,X\rangle$, so $F=\nabla f+fx$.
	Differentiating gives
	\begin{equation}\label{eq:spherical-support-differential}
		dF(X)=\nabla_X\nabla f+fX=AX.
	\end{equation}
	Conversely, for any smooth $f$, the map $F=\nabla f+fx$ has
	differential $A_f$. Its self-adjointness and $d(dF)=0$ imply Codazzi
	by the same calculation. Finally, $F$ is unique up to a constant
	vector, giving the stated freedom in $f$.
\end{proof}

We call $F=\nabla f+fx$ the \emph{support map}. The affine hyperplane
\begin{equation*}
	\{y\in\mathbb R^{m+1}:\langle y,x\rangle=f(x)\}
\end{equation*}
passes through $F(x)$ and has normal $x$, which is orthogonal to
$\operatorname{im}dF_x$. This is the support-function description of
a hedgehog, without any convexity assumption; see~\cite[Section~2]{Pan05}
for this viewpoint in $\mathbb R^3$. Although $F(\Sph^m)$ is compact,
it may have rank-drop singularities. We will use only its smooth
local images near the compact fiber components constructed below.
Also, restriction to
a totally geodesic subsphere commutes with the spherical Hessian,
so every great subsphere inherits the saddle condition.

We identify Euclidean tangent spaces with $\mathbb R^{m+1}$ by
translation. Thus the normal vectors to a local image of $F$ and the
points of the source sphere are regarded as vectors in the same space.

\begin{lemma}\label{lem:regular-fibers}
	Let $A$ be a smooth self-adjoint Codazzi tensor on $\Sph^m$, and
	let $F$ be its support map. Suppose
	$r=\max_{\Sph^m}\rank A$ satisfies $0<r<m$. There is a great
	$(m-r)$-sphere $K$ in the rank-$r$ set on which $F$ is constant.
	Along $K$, the map $F$ has a common smooth local image $S$ of
	dimension $r$. Write $F|_K\equiv z$, and put $P=T_zS$ and
	$N=P^\perp$.
	Then $K=\Sph^m\cap N$: the fiber is the entire unit normal sphere
	of $S$ at $z$.
\end{lemma}

\begin{proof}
	Choose $p\in\Sph^m$ with $\rank A(p)=r$ and a linear map
	$L:\mathbb R^{m+1}\to\mathbb R^r$ such that $L\circ A(p)$ is
	surjective. Then $G=L\circ F$ is a submersion near $p$, so its
	image contains a nonempty open subset of $\mathbb R^r$. By Sard's
	theorem, choose in this subset a regular value $y$ of the global
	map $G:\Sph^m\to\mathbb R^r$, and let $K$ be a connected component
	of $G^{-1}(y)$. It is a compact submanifold without boundary, and
	\begin{equation*}
		r=\rank dG_x\leq\rank A(x)\leq r,
		\qquad T_xK=\ker dG_x=\ker A(x)
		\quad (x\in K).
	\end{equation*}
	Thus $F|_K$ is constant; denote its value by $z$.

	The rank of $F$ is $r$ near $K$. By the constant-rank theorem,
	$F$ has an embedded local image of dimension $r$ near each point
	of $K$. If two such neighborhoods meet on $K$, their image germs
	at $z$ agree: the restriction of $F$ to their intersection still
	has rank $r$, and its image contains a neighborhood of $z$ in
	each local image.
	Connectedness of $K$ makes these image germs agree along $K$.
	A finite cover of the compact set $K$ then gives a common local
	image $S$ near $z$, after shrinking the neighborhoods.
	Consequently $P=T_zS=\operatorname{im}A(x)$ is independent of
	$x\in K$. Self-adjointness gives the ambient normal space
	\begin{equation*}
		N=\mathbb Rx\oplus\ker A(x),\qquad x\in K.
	\end{equation*}
	Since $P\subset x^\perp$, we have
	$K\subset\Sph^m\cap P^\perp$. Both sides have dimension $m-r$;
	hence $K$ is open in this great sphere, and compactness makes it
	closed. Since $m-r\geq1$, the great sphere is connected, proving
	$K=\Sph^m\cap N$.

\end{proof}

The normal-fiber description supplies a complete great circle on which
the support map is constant. To use this global information, we
restrict to a great two-sphere containing the circle and translate the
constant support value to the origin. The support function then
vanishes on the circle. This condition alone forces its restriction
to be a first spherical harmonic: Stokes' theorem makes the projected
signed area vanish on each hemisphere, and the saddle condition turns
this into equality in the first Dirichlet eigenvalue inequality.

\begin{lemma}\label{lem:circle-support-rigidity}
	Let $u\in C^2(\Sph^2)$ satisfy $\det A_u\leq0$. Suppose $u=0$
	on a complete great circle $C=\Sph^2\cap e^\perp$, where
	$e\in\mathbb R^3$ is a unit vector. Then, for some $a\in\mathbb R$,
	\begin{equation*}
		u(x)=a\langle e,x\rangle
		\qquad\text{on }\Sph^2.
	\end{equation*}
	In particular, $A_u=0$.
\end{lemma}

\begin{proof}
	Rotate the sphere so that $e=e_3$. Write $A=A_u$ and
	$F=\nabla u+ux$, and let $H_\pm=\{\pm x_3\geq0\}$ be the
	two closed hemispheres. Since $u=0$ on the equator $C$, its
	gradient there is parallel to $e_3$. Hence $F_1=F_2=0$ on $C$.
	The identity $dF=A$ gives
	\begin{equation*}
		dF_1\wedge dF_2=x_3\det A\,d\mu,
	\end{equation*}
	where $d\mu$ is the oriented round area form. This identity is
	on the source sphere and does not require $F$ to be an immersion.
	Stokes' theorem therefore yields
	\begin{equation*}
		\int_{H_\pm}x_3\det A\,d\mu
		=\int_{\partial H_\pm}F_1\,dF_2=0.
	\end{equation*}
	Since $x_3$ has a fixed nonzero sign in the interior of each
	hemisphere and $\det A\leq0$, we obtain $\det A=0$ everywhere.

	Let $H$ be either hemisphere and $\nu$ its outward unit conormal.
	The equator is totally geodesic and $u|_C=0$, so
	$\Delta u-\nabla^2u(\nu,\nu)=0$ on $C$.
	The boundary term in the integrated Bochner formula is therefore
	$u_\nu\bigl(\Delta u-\nabla^2u(\nu,\nu)\bigr)=0$, giving
	\begin{equation*}
		\int_H\bigl((\Delta u)^2-|\nabla^2u|^2\bigr)\,d\mu
		=\int_H|\nabla u|^2\,d\mu.
	\end{equation*}
	Using the two-dimensional identity
	\begin{equation*}
		2\det A=(\Delta u)^2-|\nabla^2u|^2+2u\Delta u+2u^2
	\end{equation*}
	and integrating $u\Delta u$ by parts, we conclude that
	\begin{equation*}
		0=2\int_H\det A\,d\mu
		=2\int_Hu^2\,d\mu-\int_H|\nabla u|^2\,d\mu.
	\end{equation*}
	The integral identities above are valid for $C^2$ functions by
	smooth approximation; their interior and boundary terms involve
	derivatives of order at most two.

	The first Dirichlet eigenvalue of a unit hemisphere is $2$, and
	its eigenspace is spanned by $x_3$. The last identity is precisely
	the equality case in the Rayleigh inequality. Thus
	$u=a_+x_3$ on $H_+$ and $u=a_-x_3$ on $H_-$.
	Continuity of $\nabla u$ across $C$, where $\nabla x_3=e_3$,
	gives $a_+=a_-$. This proves the claim.
\end{proof}

\begin{theorem}\label{prop:spherical-saddle-rigidity}
	Let $m\geq3$. Every smooth self-adjoint saddle Codazzi tensor on $\Sph^m$
	 vanishes identically. 	Equivalently, if $f\in C^\infty(\Sph^m)$ and $A_f$ is saddle, then	$f(x)=\langle z,x\rangle$ for some $z\in\mathbb R^{m+1}$.
\end{theorem}

\begin{proof}
	Write $A=A_f$ and $F=\nabla f+fx$. Suppose $A\not\equiv0$.
	By Lemma~\ref{lem:rank-two}, its maximal rank satisfies
	$0<r\leq2<m$. Lemma~\ref{lem:regular-fibers} therefore gives a
	great $(m-r)$-sphere $K$ on which $F$ is constant, say $F=z$.
	Since $m-r\geq1$, choose a complete great circle $C\subset K$.
	For $\widehat f(x)=f(x)-\langle z,x\rangle$, we have
	\begin{equation*}
		\widehat f=0,\qquad d\widehat f=0\qquad\text{on }C.
	\end{equation*}

	Let $\Sigma\simeq\Sph^2$ be any great two-sphere containing $C$,
	and set $u=\widehat f|_\Sigma$. Since $\Sigma$ is totally geodesic,
	$A_u$ represents the restriction of the bilinear form associated
	with $A$ to $T\Sigma$. Hence $\det A_u\leq0$, and $u=0$ on $C$.
	Lemma~\ref{lem:circle-support-rigidity} makes $u$ a linear
	function on $\Sigma$. Since also $du=0$ on $C$, this linear
	function is zero. Thus $u=0$ on $\Sigma$.
	Every point of $\Sph^m$ lies on a great two-sphere containing $C$:
	adjoin its position vector to the two-plane spanning $C$, and
	extend to a three-plane if necessary. Thus $\widehat f=0$ on
	$\Sph^m$, so $A=0$, a contradiction.
\end{proof}

The dimension threshold is now explicit. For $m\geq3$, the bound
$r\leq2<m$ supplies a positive-dimensional normal fiber and hence a
complete circle. On $\Sph^2$, maximal rank two supplies no such circle;
this is the case in which the nontrivial saddle data used below occur.

\begin{proof}[Proof of Theorem~\ref{thm:main}]
	Identify $\partial M$ with $\Sph^{n-1}$. Its second fundamental
	form is Codazzi by assumption, and
	\eqref{eq:boundary-saddle-gauss} gives the saddle condition for
	the associated self-adjoint endomorphism.
	Theorem~\ref{prop:spherical-saddle-rigidity} therefore
	gives $h=0$. Since $\sec_g\geq1$ implies $\Ric_g\geq(n-1)g$,
	the hemisphere rigidity theorem of Hang--Wang~\cite{HW} gives
	$(M^n,g)\cong\Sph^n_+$.
\end{proof}

\section{Construction of counterexamples in dimension three}\label{sec:counterexamples}

Unlike the case when $m\geq 3$, there do exist nontrivial saddle Codazzi tensors on $\Sph^2$, see e.g., ~\cite[Theorem~6.1]{Pan05}.
We shall realize the nontrivial saddle data available on $\Sph^2$,
proving Theorem~\ref{thm:three}.

\subsection{Saddle data and the metric realization}

 Let $H=\Sph^3_+$ be the closed northern hemisphere, equipped with
the standard round metric $g_0$, and let $\sigma$ be the induced
unit round metric on $\partial H=\Sph^2$.  Panina's construction provides a function
$f\in C^\infty(\Sph^2)$ such that
\begin{equation}\label{eq:counterexample-saddle-data}
	A=\nabla_\sigma^2 f+f\sigma,
	\qquad
	\rho=-\det_\sigma A\geq 0,
	\qquad
	\rho\not\equiv 0.
\end{equation}
We also use $A$ for its associated self-adjoint endomorphism with
respect to $\sigma$.
The associated support map is
\begin{equation*}
	F_f(q)=\nabla_\sigma f(q)+f(q)q,
	\qquad dF_f=A.
\end{equation*}
The condition $\rho\not\equiv0$ follows already from
$A\not\equiv0$. Indeed, if $\rho\equiv0$, then $A$ has maximal
rank one. Lemma~\ref{lem:regular-fibers}, with $m=2$ and $r=1$,
gives a complete great circle on which $F_f$ is constant, say $z$.
The function $f-\langle z,q\rangle$ vanishes on that circle, so
Lemma~\ref{lem:circle-support-rigidity} would give $A=0$.
Panina's construction gives a smooth function and hence a smooth map
$F_f$; the map need not be an immersion, and its image need not be
a smooth surface. We use $F_f$ to organize auxiliary Euclidean data.
The metric on $H$ will be produced by a collar construction and an
interior correction that control its sectional curvature globally.
In particular, no realization as the boundary of a Euclidean convex
body is assumed. The tensor $A$ is Codazzi by
Lemma~\ref{lem:codazzi-potential-map}.
For the proposed boundary data $h_\varepsilon=\varepsilon A$, the
tangential Gauss equation requires
\begin{equation*}
	\sec_{g_\varepsilon}(T\partial H)=1-\varepsilon^2\det_\sigma A
	=1+\varepsilon^2\rho\geq1.
\end{equation*}
The task is to realize these data while controlling all sectional
curvatures in the interior.



For a metric $g$ on $H$, we use the convention
\begin{equation*}
	h_g(X,Y)=g(\nabla^g_X\nu_g,Y),
	\qquad X,Y\in T\partial H,
\end{equation*}
where $\nu_g$ is the outward unit normal.

We shall prove the following.

\begin{theorem}\label{thm:three-dimensional-counterexample}
	For every $f\in C^\infty(\Sph^2)$ satisfying
	\eqref{eq:counterexample-saddle-data}, there exist $\varepsilon_0>0$
	and a smooth family of Riemannian metrics
	$\{g_\varepsilon\}_{|\varepsilon|<\varepsilon_0}$ on $H$ such that
	\begin{equation*}
		g_{\varepsilon=0}=g_0,
		\qquad
		\sec_{g_\varepsilon}\geq 1,
		\qquad
		g_\varepsilon|_{T\partial H}=\sigma,
		\qquad
		h_{g_\varepsilon}=\varepsilon A.
	\end{equation*}
	In particular, the boundary second fundamental form is Codazzi for
	every $|\varepsilon|<\varepsilon_0$, and $(H,g_\varepsilon)$ is not
	isometric to the standard hemisphere whenever
	$0<|\varepsilon|<\varepsilon_0$.
\end{theorem}

Away from the north pole, we write
\begin{equation*}
	X=(x,z)=(\cos t\,q,\sin t)\in \Sph^3\subset\mathbb R^3\times\mathbb R,
	\qquad
	g_0=dt^2+\cos^2t\,\sigma,
\end{equation*}
where $q\in\Sph^2$ and $t\in[0,\pi/2)$ is the inward distance from
the equator. Near the equator, we also use the signed extension of
$t$ to the ambient sphere, with $t<0$ on the southern side.
On the open northern hemisphere, the gnomonic coordinates are
$
y=\frac{x}{z}\in\mathbb R^3.
$
For $0<t<\pi/2$, this gives $y=\cot t\,q$, while the north pole
corresponds to $y=0$. In these coordinates,
\begin{equation*}
	g_0
	=
	\frac{\sum_i(dy^i)^2}{1+|y|^2}
	-
	\frac{\bigl(\sum_i y^i\,dy^i\bigr)^2}{(1+|y|^2)^2}.
\end{equation*}
Thus, for $0<c<\pi/2$, the interior cap $\{t\geq c\}$ corresponds
to the closed Euclidean ball $\{|y|\leq\cot c\}$.

In dimension three, the sectional curvature bound can be expressed
as the nonnegativity of a symmetric two-tensor. For a Riemannian
metric $g$, define
\begin{equation*}
	T(g)
	=
	\left(\frac{\operatorname{Scal}_g}{2}-1\right)g
	-\operatorname{Ric}_g.
\end{equation*}
If $n$ is a $g$-unit vector, then
\begin{equation*}
	T(g)(n,n)=\sec_g(n^\perp)-1,
\end{equation*}
where $n^\perp$ is the two-plane orthogonal to $n$ with respect to
$g$. Consequently,
\begin{equation*}
	\sec_g\geq1
	\quad\Longleftrightarrow\quad
	T(g)\geq0.
\end{equation*}
More generally, $T(g)\geq\delta g$ is equivalent to
$\sec_g\geq1+\delta$. Thus $T(g)$ records the sectional
curvature excess over $1$. Moreover, the contracted Bianchi identity
gives
\begin{equation*}
	\operatorname{div}_gT(g)=0,
	\qquad
	T(g_0)=0.
\end{equation*}

The construction has three steps. Choose the parameters in the order
described below, and then fix them independently of $\varepsilon$:

\begin{equation*}
	0<t_0<a<b<t_1<t_2<\frac{\pi}{4},
\end{equation*}
First, using gnomonic coordinates, we construct a smooth symmetric
tensor $D$ on an open neighborhood of the closed hemisphere $H$
such that
\begin{equation*}
	D\geq0,
	\qquad
	\operatorname{div}_{g_0}D=0.
\end{equation*}
This tensor is intended to be the coefficient of $\varepsilon^2$
in $T(g_\varepsilon)$. Its normal components at the boundary
are prescribed by the Gauss and Codazzi equations:
\begin{equation*}
	D(\partial_t,\partial_t)=\rho,
	\qquad
	D(\partial_t,X)=0
	\quad\text{for }X\in T\partial H.
\end{equation*}
We also arrange that $D$ is uniformly positive definite on the
interior cap $\{t\geq t_0\}$. On a larger, two-sided collar
$[-\tau,t_2]\times\Sph^2$, where $\tau>0$, the tensor $D$ admits
a nonnegative geodesic integral representation, which provides
the data for the next step.

Second, we construct a smooth family of metrics $G_\varepsilon$
on the fixed boundary collar $\{0\leq t\leq t_1\}$, with induced
boundary metric $\sigma$ and second fundamental form
$\varepsilon A$, such that
\begin{equation*}
	T(G_\varepsilon)\geq0,
	\qquad
	T(G_\varepsilon)
	=
	\varepsilon^2D+O(\varepsilon^3).
\end{equation*}
The first inequality is exact, including at points where $D$ is
degenerate. On the interior slab $\{t_0\leq t\leq t_1\}$, the
positivity of $D$ gives a strictly positive curvature margin of
order $\varepsilon^2$.

Finally, we extend the collar metrics smoothly to $H$, retaining
$G_\varepsilon$ on $\{t\leq a\}$ and placing the cutoff transition
inside the slab $\{a<t<b\}$. The extension is arranged so that its
first variation is a Lie derivative of $g_0$. Since
$T(g_0)=0$ and $T$ is natural under
diffeomorphisms, its curvature excess has no linear term in
$\varepsilon$. The discrepancy between its quadratic term and
$D$ is divergence-free and supported in the closed cap
$\{t\geq a\}$. A supported solution of the linearized equation
for $T$ removes this discrepancy by a metric correction
of order $\varepsilon^2$, giving
\begin{equation*}
	T(g_\varepsilon)
	=
	\varepsilon^2D+O(\varepsilon^3)
	\quad\text{on }H.
\end{equation*}
On $\{t\geq a\}$, the uniform positivity of $D$ absorbs the
remaining cubic error for sufficiently small $|\varepsilon|$.
On $\{t\leq a\}$, the metric remains equal to the exact collar
metric, so the curvature bound and the prescribed boundary data
are preserved. All remainder estimates are understood in fixed
$C^k$ norms on the indicated regions.

\subsection{Curvature data and geodesic transport}

We first extend the boundary curvature excess to a nonnegative
divergence-free tensor. Its Euclidean description gives both the
required positivity in the interior and a representation by round
geodesics, which will persist under the collar deformation.

\begin{lemma}\label{lem:divergence-free-tensor}
	There exists a smooth symmetric two-tensor $D$ on $H$, smooth up
	to the boundary, satisfying
	\begin{equation*}
		D\geq0,
		\qquad
		\operatorname{div}_{g_0}D=0.
	\end{equation*}
	At the equator,
	\begin{equation*}
		D(\partial_t,\partial_t)=\rho,
		\qquad
		D(\partial_t,Y)=0
		\quad\text{for every }Y\in T\partial H.
	\end{equation*}
	Moreover, there exist constants $c_0>0$ and $0<t_0<\pi/4$ such that
	\begin{equation*}
		D\geq c_0g_0
		\qquad\text{on }\{t\geq t_0\}.
	\end{equation*}
\end{lemma}

\begin{proof}
	We first establish the balance identity
	\begin{equation}\label{eq:saddle-balance}
		\int_{\Sph^2}\rho(q)\,q\,d\mu_\sigma(q)=0.
	\end{equation}
	The identity states that the total oriented area vector of $F_f$
	vanishes. We can apply Stokes' theorem on the parameter sphere
	even where $F_f$ is singular. Orient $\Sph^2$ by its outward normal $q$.
	For every cyclic permutation $(i,j,k)$ of $(1,2,3)$, the identity
	$dF_f=A$ gives
	\begin{equation*}
		dF_f^j\wedge dF_f^k
		=(\det_\sigma A)q^i\,d\mu_\sigma.
	\end{equation*}
	Indeed, on a positively oriented orthonormal frame, the cross
	product of the two images under $A$ is $(\det_\sigma A)q$.
	Stokes' theorem now yields
	\begin{equation*}
		\int_{\Sph^2}(\det_\sigma A)q^i\,d\mu_\sigma
		=\int_{\Sph^2}dF_f^j\wedge dF_f^k
		=\int_{\Sph^2}d\bigl(F_f^j\,dF_f^k\bigr)=0.
	\end{equation*}
	Since $\rho=-\det_\sigma A$, this proves
	\eqref{eq:saddle-balance}.

	We next construct a tensor in Euclidean space. For $y\neq0$, define
	the symmetric contravariant tensor
	\begin{equation}\label{eq:euclidean-seed-tensor}
		Q^{ij}(y)
		=
		\rho\!\left(\frac{y}{|y|}\right)
		\frac{y^iy^j}{|y|^4}.
	\end{equation}
	It is nonnegative, locally integrable, and satisfies
	$\partial_iQ^{ij}=0$ away from the origin. Its distributional
	divergence also vanishes at the origin. Indeed, for a compactly
	supported smooth vector field $\varphi$ on $\mathbb R^3$, integration
	by parts outside a ball of radius $r$ gives
	\begin{equation*}
		\begin{aligned}
		\langle\operatorname{div}Q,\varphi\rangle
		&=
		\lim_{r\downarrow0}
		\int_{\Sph^2}
		\rho(q)\,\langle q,\varphi(rq)\rangle\,d\mu_\sigma(q)\\
		&=
		\left\langle
		\int_{\Sph^2}\rho(q)q\,d\mu_\sigma(q),\varphi(0)
		\right\rangle
		=0.
		\end{aligned}
	\end{equation*}

	Fix a smooth radial function $\eta\geq0$ on $\mathbb R^3$ such that
	\begin{equation*}
		\operatorname{supp}\eta\subset\overline{B_1(0)},
		\qquad
		\eta>0\ \text{on }B_1(0),
		\qquad
		\int_{\mathbb R^3}\eta(a)\,da=1.
	\end{equation*}
	The convolution
	\begin{equation*}
		\widehat D^{ij}(y)
		=
		(Q^{ij}*\eta)(y)
		=
		\int_{\mathbb R^3}
		\rho\!\left(\frac{y-a}{|y-a|}\right)
		\frac{(y-a)^i(y-a)^j}{|y-a|^4}\,\eta(a)\,da
	\end{equation*}
	is smooth, nonnegative, and divergence-free with respect to the
	Euclidean metric.

	We now transfer this tensor to the open northern hemisphere using
	the gnomonic coordinates from the setup. Write $s=1+|y|^2$, and
	define $D$ by prescribing its contravariant components:
	\begin{equation}\label{eq:gnomonic-tensor-definition}
		D^{ij}=s^3\widehat D^{ij}.
	\end{equation}
	Here the indices of $D$ are raised and lowered using $g_0$.
	In particular, its covariant components are
	\begin{equation*}
		D_{ij}
		=
		s^3(g_0)_{ik}(g_0)_{j\ell}\widehat D^{k\ell},
	\end{equation*}
	so $D\geq0$. Writing $y_i=y^i$ for the Euclidean coordinate
	components, the Christoffel symbols are
	\begin{equation*}
		\Gamma^k_{ij}
		=
		-\frac{y_i\delta^k_j+y_j\delta^k_i}{s}.
	\end{equation*}
	Thus, for a symmetric contravariant tensor,
	\begin{equation*}
		\nabla_iD^{ij}
		=
		\partial_iD^{ij}-\frac{6y_i}{s}D^{ij}
		=
		s^3\partial_i\bigl(s^{-3}D^{ij}\bigr).
	\end{equation*}
	Together with \eqref{eq:gnomonic-tensor-definition}, this proves
	$\operatorname{div}_{g_0}D=0$ in the interior of $H$.

	To verify smoothness at the equator, we rewrite the construction
	in ambient spherical coordinates. For $a\in\mathbb R^3$, set
	\begin{equation*}
		w_a=x-za,
		\qquad
		W_a=\operatorname{pr}_{T_X\Sph^3}(a,1).
	\end{equation*}
	The tensor defined above has the equivalent expression
	\begin{equation}\label{eq:D-spherical-representation}
		D_X
		=
		\int_{\mathbb R^3}
		\rho\!\left(\frac{w_a}{|w_a|}\right)
		\frac{W_a^\flat\otimes W_a^\flat}{|w_a|^4}\,
		\eta(a)\,da,
	\end{equation}
	where $\flat$ denotes the dual with respect to $g_0$.
	Indeed, in the northern gnomonic coordinates,
	\begin{equation*}
		w_a=\frac{y-a}{\sqrt{s}},
		\qquad
		W_a^i=\sqrt{s}\,(a^i-y^i),
	\end{equation*}
	which gives precisely \eqref{eq:gnomonic-tensor-definition}.

	For $a\in\operatorname{supp}\eta$ and $|t|<\pi/4$,
	\begin{equation*}
		|w_a|
		\geq \cos t-|\sin t|>0.
	\end{equation*}
	Hence \eqref{eq:D-spherical-representation} gives a smooth,
	nonnegative extension of $D$ across the equator. In particular,
	$D$ is smooth up to $\partial H$, and its divergence vanishes there
	by continuity.

	At $X=(q,0)$, we have
	\begin{equation*}
		w_a=q,
		\qquad
		W_a=\bigl(a-\langle a,q\rangle q,1\bigr),
		\qquad
		\partial_t=(0,1).
	\end{equation*}
	It follows that
	\begin{equation*}
		D(\partial_t,\partial_t)
		=
		\rho(q)\int_{\mathbb R^3}\eta(a)\,da
		=
		\rho(q).
	\end{equation*}
	For $Y\in T_q\Sph^2$, radial symmetry of $\eta$ gives
	\begin{equation*}
		D(\partial_t,Y)
		=
		\rho(q)\int_{\mathbb R^3}\langle a,Y\rangle\eta(a)\,da
		=
		0.
	\end{equation*}

	It remains to prove strict positivity on an interior cap.
	We claim that $\widehat D(y)$ is positive definite for every
	$|y|\leq1$. Fix such a point $y$ and a nonzero Euclidean covector
	$\xi$. If $|y|<1$, choose $q\in\Sph^2$ such that
	\begin{equation*}
		\rho(q)>0,
		\qquad
		\langle\xi,q\rangle\neq0.
	\end{equation*}
	Such a choice is possible because $\{\rho>0\}$ is a nonempty open
	set. For sufficiently small $r>0$, the point $a=y-rq$ lies in
	$B_1(0)$.

	If $|y|=1$, the balance identity \eqref{eq:saddle-balance} implies
	that the open set
	\begin{equation*}
		\{q\in\Sph^2:\rho(q)>0,\ \langle y,q\rangle>0\}
	\end{equation*}
	is nonempty. Indeed, otherwise the nonnegative smooth density
	$\rho$ would be supported in a closed hemisphere, contradicting
	\eqref{eq:saddle-balance}. We may therefore choose $q$ in this set
	with $\langle\xi,q\rangle\neq0$. Again, $a=y-rq$ belongs to $B_1(0)$
	for sufficiently small $r>0$.

	In either case, the integrand in
	\begin{equation*}
		\widehat D^{ij}(y)\xi_i\xi_j
		=
		\int_{\mathbb R^3}
		\rho\!\left(\frac{y-a}{|y-a|}\right)
		\frac{\langle\xi,y-a\rangle^2}{|y-a|^4}\,
		\eta(a)\,da
	\end{equation*}
	is strictly positive on a nonempty open set. Thus
	$\widehat D(y)$ is positive definite on $\overline{B_1(0)}$.
	By smoothness and compactness, it remains positive definite on
	$\overline{B_R(0)}$ for some $R>1$. Since $|y|=\cot t$, this ball
	corresponds to the cap $\{t\geq t_0\}$, where
	$
	t_0=\operatorname{arccot}R<\frac{\pi}{4}.
	$
	The asserted uniform lower bound follows from compactness of this
	cap.
\end{proof}

We next show that the extension
\eqref{eq:D-spherical-representation} has a representation adapted
to the round geodesic flow. This additional structure will be used
in the construction of the collar metric.

\begin{lemma}\label{lem:geodesic-representation}
	Let $D$ be the tensor constructed in
	Lemma~\ref{lem:divergence-free-tensor}, extended across the
	equator by \eqref{eq:D-spherical-representation}.
	Choose $t_0<t_1<t_2<\pi/4$. There exist $\tau>0$, $\kappa>0$,
	and a smooth nonnegative function $F_{\mathrm{rd}}$ on the round
	unit tangent bundle over a neighborhood of
	\begin{equation*}
		C_c=[-\tau,t_2]\times\Sph^2
	\end{equation*}
	such that
	\begin{equation}\label{eq:D-geodesic-representation}
		D_X
		=
		\int_{S_X\Sph^3}
		F_{\mathrm{rd}}(X,v)\,
		v^\flat\otimes v^\flat\,d\mu_X(v).
	\end{equation}
	The function $F_{\mathrm{rd}}$ is invariant under the round
	geodesic flow, as long as the base geodesic remains in this
	neighborhood, and
	\begin{equation*}
		\operatorname{supp}F_{\mathrm{rd}}\big|_{C_c}
		\subset
		\{(X,v):v^0\geq\kappa\},
		\qquad
		v^0=\langle v,\partial_t\rangle_{g_0}.
	\end{equation*}
	At the equator,
	\begin{equation}\label{eq:geodesic-weight-boundary}
		F_{\mathrm{rd}}(0,q,v)=\rho(q)\Phi(v),
	\end{equation}
	where $\Phi$ is smooth, nonnegative, invariant under rotations
	of the tangential components, and satisfies
	\begin{equation}\label{eq:geodesic-weight-moments}
		\int_{\Sph^2}\Phi(v)(v^0)^2\,d\mu(v)=1,
		\qquad
		\int_{\Sph^2}\Phi(v)v^0v^\alpha\,d\mu(v)=0,
		\quad \alpha=1,2.
	\end{equation}
	Here $v^\alpha$ are the components in any orthonormal tangential
	frame, and the integrals are over the sphere of unit tangent
	directions.
\end{lemma}

\begin{proof}
	We work on an open collar $U$ containing $C_c$, with
	$\overline U\subset\{|t|<\pi/4\}$. For
	$a\in\operatorname{supp}\eta$, retain the notation $w_a,W_a$ from
	\eqref{eq:D-spherical-representation}, and set
	\begin{equation*}
		\ell_a(X)=\langle X,(a,1)\rangle,
		\qquad
		r_a=|W_a|_{g_0},
		\qquad
		v_a=\frac{W_a}{r_a},
	\end{equation*}
	and
	\begin{equation*}
		\lambda_a
		=
		\rho\!\left(\frac{w_a}{|w_a|}\right)
		\frac{r_a^2}{|w_a|^4}.
	\end{equation*}
	These quantities are smooth on $U$. Indeed,
	\begin{equation*}
		|w_a|\geq\cos t-|\sin t|,
		\qquad
		\langle W_a,\partial_t\rangle
		=
		\cos t-\sin t\,\langle q,a\rangle
		\geq\cos t-|\sin t|.
	\end{equation*}
	In particular, $w_a$ and $W_a$ are nonzero. Since
	$r_a\leq|(a,1)|\leq\sqrt2$, we also obtain a uniform lower bound
	\begin{equation*}
		\langle v_a,\partial_t\rangle\geq\kappa>0
		\qquad\text{on }C_c.
	\end{equation*}

	We first verify the geometric properties of each vector field
	$v_a$. Since $W_a=\nabla^{g_0}\ell_a$ and $\ell_a$ is the
	restriction of a linear function to the unit sphere,
	\begin{equation*}
		\nabla_ZW_a=-\ell_a Z
		\qquad\text{for every }Z\in T\Sph^3.
	\end{equation*}
	Direct differentiation also gives
	\begin{equation*}
		W_a(w_a)=-\ell_a w_a,
		\qquad
		W_a(r_a)=-\ell_a r_a.
	\end{equation*}
	It follows that
	\begin{equation*}
		\nabla_{v_a}v_a=0,
		\qquad
		\operatorname{div}_{g_0}v_a=-\frac{2\ell_a}{r_a},
		\qquad
		v_a(\lambda_a)=\frac{2\ell_a}{r_a}\lambda_a.
	\end{equation*}
	Thus the integral curves of $v_a$ are unit-speed round geodesics,
	and their weighted flux satisfies
	\begin{equation}\label{eq:geodesic-beam-conservation}
		\operatorname{div}_{g_0}(\lambda_av_a)=0.
	\end{equation}
	These identities also imply
	\begin{equation*}
		\operatorname{div}_{g_0}
		\bigl(\lambda_a v_a^\flat\otimes v_a^\flat\bigr)=0.
	\end{equation*}
	Consequently, the extension of $D$ is divergence-free on the
	whole two-sided collar.

	For $X\in U$, push the measure $\lambda_a(X)\eta(a)\,da$ forward
	by $a\mapsto v_a(X)$ to the unit tangent sphere. Denote the
	resulting density by $F_{\mathrm{rd}}$:
	\begin{equation}\label{eq:geodesic-density-pushforward}
		\int_{S_X\Sph^3}\psi(v)F_{\mathrm{rd}}(X,v)\,d\mu_X(v)
		=
		\int_{\mathbb R^3}
		\psi(v_a(X))\lambda_a(X)\eta(a)\,da.
	\end{equation}
	To justify that this density is smooth, write
	\begin{equation*}
		a=u+bq,
		\qquad
		u\in T_q\Sph^2.
	\end{equation*}
	Relative to the normal and tangential splitting at $X$, the
	components of $W_a$ are
	\begin{equation*}
		W_a=(\cos t-b\sin t,u).
	\end{equation*}
	On a neighborhood of $\operatorname{supp}\eta$, its first component
	is positive. For fixed $b$, the map
	\begin{equation*}
		u\longmapsto
		\frac{(\cos t-b\sin t,u)}
		{\sqrt{(\cos t-b\sin t)^2+|u|^2}}
	\end{equation*}
	has rank two. Hence $a\mapsto v_a(X)$ is a submersion on a
	neighborhood of the support of the integrand. Integration along
	its fibers therefore gives a smooth density, depending smoothly
	on $X$. This density is nonnegative, and its support satisfies
	the stated transversality condition.
	Moreover, since
	\begin{equation*}
		\lambda_a v_a^\flat\otimes v_a^\flat
		=
		\rho\!\left(\frac{w_a}{|w_a|}\right)
		\frac{W_a^\flat\otimes W_a^\flat}{|w_a|^4},
	\end{equation*}
	formula \eqref{eq:geodesic-density-pushforward} gives
	\eqref{eq:D-geodesic-representation}.

	We now prove invariance under the round geodesic flow. Let
	$X_0$ be its generating vector field on the unit tangent
	bundle $SU$. For every $\psi\in C_c^\infty(SU)$, the identity
	$\nabla_{v_a}v_a=0$ gives
	\begin{equation*}
		(X_0\psi)(X,v_a(X))
		=
		v_a\bigl[\psi(X,v_a(X))\bigr].
	\end{equation*}
	Using \eqref{eq:geodesic-density-pushforward} and
	\eqref{eq:geodesic-beam-conservation}, we obtain
	\begin{equation*}
		\begin{aligned}
		&\int_{SU}
		F_{\mathrm{rd}}\,X_0\psi\,
		d\operatorname{vol}_{g_0}(X)\,d\mu_X(v)\\
		&\qquad=
		\int_{\mathbb R^3}\eta(a)
		\int_U
		\lambda_a\,v_a\bigl[\psi(X,v_a(X))\bigr]\,
		d\operatorname{vol}_{g_0}(X)\,da\\
		&\qquad=0.
		\end{aligned}
	\end{equation*}
	The round geodesic flow preserves the Liouville measure
	$d\operatorname{vol}_{g_0}(X)\,d\mu_X(v)$. Integration by parts
	therefore yields $
	X_0F_{\mathrm{rd}}=0.
	$
	Notice that this invariance concerns the density
	$F_{\mathrm{rd}}$ on the unit tangent bundle; the individual
	weights $\lambda_a$ satisfy the flux conservation law
	\eqref{eq:geodesic-beam-conservation}.

	Finally, we compute the density at the equator.
	Writing again $a=u+bq$, we have, in normal and tangential
	components,
	\begin{equation*}
		W_a=(1,u),
		\qquad
		v_a=\frac{(1,u)}{\sqrt{1+|u|^2}},
		\qquad
		\lambda_a=\rho(q)(1+|u|^2).
	\end{equation*}
	Define
	\begin{equation*}
		\eta_{\mathrm{tan}}(u)
		=
		\int_{\mathbb R}\eta(u+bq)\,db.
	\end{equation*}
	Radial symmetry of $\eta$ implies that $\eta_{\mathrm{tan}}$
	depends only on $|u|$ and is independent of $q$ after an
	orthonormal identification of $T_q\Sph^2$ with $\mathbb R^2$.
	Under the parametrization
	\begin{equation*}
		v=\frac{(1,u)}{\sqrt{1+|u|^2}},
	\end{equation*}
	the area element on the unit tangent sphere is
	\begin{equation*}
		d\mu(v)=(1+|u|^2)^{-3/2}\,du.
	\end{equation*}
	Consequently, \eqref{eq:geodesic-weight-boundary} holds with
	\begin{equation*}
		\Phi\!\left(\frac{(1,u)}{\sqrt{1+|u|^2}}\right)
		=
		(1+|u|^2)^{5/2}\eta_{\mathrm{tan}}(u).
	\end{equation*}
	Since $\eta_{\mathrm{tan}}$ is smooth and supported in
	$\{|u|\leq1\}$, this formula extends by zero to a smooth function
	on the entire unit sphere, supported where $v^0\geq1/\sqrt2$.
	Its rotational invariance follows from that of
	$\eta_{\mathrm{tan}}$. Finally,
	\begin{equation*}
		\int_{\Sph^2}\Phi(v)(v^0)^2\,d\mu(v)
		=
		\int_{\mathbb R^2}\eta_{\mathrm{tan}}(u)\,du
		=
		1,
	\end{equation*}
	and
	\begin{equation*}
		\int_{\Sph^2}\Phi(v)v^0v^\alpha\,d\mu(v)
		=
		\int_{\mathbb R^2}u^\alpha\eta_{\mathrm{tan}}(u)\,du
		=
		0.
	\end{equation*}
	This proves \eqref{eq:geodesic-weight-moments}.
\end{proof}

\subsection{Realization on an exact boundary collar}

We now realize the tensor $D$ as the leading curvature excess of
a collar metric with the prescribed boundary data. The geodesic
representation from Lemma~\ref{lem:geodesic-representation} is
essential near the equator: transporting its nonnegative density
will give an exact curvature inequality even where $D$ vanishes.
We write $\operatorname{Lie}_Vg$ for the Lie derivative of a metric
$g$ along a vector field $V$.

\begin{lemma}\label{lem:exact-collar}
	There exist $\varepsilon_0>0$ and a fixed open collar
	$
	U\supset[0,t_1]\times\Sph^2
	$
	such that, for $|\varepsilon|<\varepsilon_0$, there is a smooth
	family of metrics
	$
	G_\varepsilon=dt^2+\gamma_\varepsilon(t)
	$
	on $U$ satisfying
	\begin{equation*}
		G_{\varepsilon=0}=g_0,
		\qquad
		\sec_{G_\varepsilon}\geq1,
		\qquad
		\gamma_\varepsilon(0)=\sigma,
		\qquad
		\dot\gamma_\varepsilon(0)=-2\varepsilon A.
	\end{equation*}
	Moreover,
	\begin{equation}\label{eq:collar-expansions}
		T(G_\varepsilon)
		=
		\varepsilon^2D+O(|\varepsilon|^3),
		\qquad
		G_\varepsilon
		=
		g_0+\varepsilon\operatorname{Lie}_Vg_0+O(\varepsilon^2),
	\end{equation}
	where
	$
	V=f\partial_t-\tan t\,\nabla_\sigma f.
	$
	In particular, there exists $c>0$ such that
	\begin{equation*}
		\sec_{G_\varepsilon}\geq1+c\varepsilon^2
		\qquad\text{on }[t_0,t_1]\times\Sph^2
	\end{equation*}
	for all sufficiently small $|\varepsilon|$.
	The expansions hold in every fixed $C^k$ norm on the collar.
\end{lemma}

\begin{proof}
	We construct the metric and an auxiliary metric connection
	simultaneously. The evolution equations will impose the curvature
	components involving the normal direction. We will then verify,
	using a conservation law and the Bianchi identities, that the
	remaining curvature equations also hold and that the auxiliary
	connection is the Levi--Civita connection.

	\medskip
	\noindent
	\emph{Step 1: the unknown metric and the evolution equations.}\par\nopagebreak[4]

	Choose a local oriented $\sigma$-orthonormal coframe
	$\{\omega_\sigma^1,\omega_\sigma^2\}$, pulled back to the collar.
	We seek spatial one-forms
	\begin{equation*}
		\omega^\alpha
		=
		(C_\varepsilon)^\alpha{}_\beta\omega_\sigma^\beta,
		\qquad \alpha,\beta\in\{1,2\},
	\end{equation*}
	and set
	\begin{equation*}
		\omega^0=dt,
		\qquad
		G_\varepsilon
		=
		\sum_{a=0}^2\omega^a\otimes\omega^a
		=
		dt^2+\sigma(C_\varepsilon\cdot,C_\varepsilon\cdot).
	\end{equation*}
	Thus the matrix $C_\varepsilon$ determines the metric on each
	slice. It suffices that $C_\varepsilon$ remain invertible.
	All time derivatives of spatial forms are taken after pulling
	them back to the fixed sphere $\Sph^2$.

	Let $\{e_a\}_{a=0}^2$ be the dual frame, so that
	$e_0=\partial_t$. We also seek a nonnegative function
	$F=F_\varepsilon(t,q,v)$ on unit tangent directions. Writing
	\begin{equation*}
		v=v^ae_a,
		\qquad
		\sum_{a=0}^2(v^a)^2=1,
	\end{equation*}
	define
	\begin{equation}\label{eq:collar-moment-matrix}
		M_{ab}(t,q)
		=
		\int_{\Sph^2_v}
		F(t,q,v)\,v^av^b\,d\mu(v).
	\end{equation}
	Here $\Sph^2_v$ is the sphere of velocity components, and the
	integration takes place at the fixed base point $(t,q)$.
	The matrix $M$ is symmetric and nonnegative, since
	\begin{equation*}
		M_{ab}\xi^a\xi^b
		=
		\int_{\Sph^2_v}
		F(t,q,v)
		\left(\sum_a\xi^av^a\right)^2\,d\mu(v)
		\geq0.
	\end{equation*}

	Our objective is to construct $G_\varepsilon$ so that
	\begin{equation}\label{eq:collar-exact-curvature-target}
		T(G_\varepsilon)
		=
		\varepsilon^2M_{ab}\,\omega^a\otimes\omega^b.
	\end{equation}
	Indeed, for a $G_\varepsilon$-unit vector $n=n^ae_a$, this identity
	would give
	\begin{equation*}
		\sec_{G_\varepsilon}(n^\perp)-1
		=
		\varepsilon^2
		\int_{\Sph^2_v}
		F(t,q,v)
		\left(\sum_an^av^a\right)^2\,d\mu(v)
		\geq0.
	\end{equation*}

	To obtain a closed evolution system, we introduce an auxiliary
	metric connection with connection forms $\Gamma^a{}_b$.
	At this stage it is not assumed to be torsion-free.
	We impose
	\begin{equation*}
		\Gamma^\alpha{}_0=\dot\omega^\alpha,
		\qquad
		\Gamma^a{}_b=-\Gamma^b{}_a,
		\qquad
		\Gamma^a{}_b(\partial_t)=0,
	\end{equation*}
	and regard the remaining spatial form $\Gamma^1{}_2$ as an
	independent unknown. We use the coefficient convention
	\begin{equation*}
		\Gamma^b{}_c=\Gamma^b{}_{ca}\omega^a,
		\qquad
		\nabla^\Gamma_{e_a}e_c=\Gamma^b{}_{ca}e_b.
	\end{equation*}

	The curvature forms corresponding to the desired curvature
	operator are prescribed by
	\begin{equation}\label{eq:collar-prescribed-curvature}
		\begin{pmatrix}
		K^1{}_2\\
		K^2{}_0\\
		K^0{}_1
		\end{pmatrix}
		=
		(I+\varepsilon^2M)
		\begin{pmatrix}
		\omega^1\wedge\omega^2\\
		\omega^2\wedge\omega^0\\
		\omega^0\wedge\omega^1
		\end{pmatrix}.
	\end{equation}
	The three rows correspond to the planes with normals
	$e_0,e_1,e_2$, respectively. Once these are the actual curvature
	forms of $G_\varepsilon$, this prescription is equivalent to
	\eqref{eq:collar-exact-curvature-target}.

	First require that
	\begin{equation*}
		d\Gamma^a{}_b+\Gamma^a{}_c\wedge\Gamma^c{}_b
	\end{equation*}
	agree with $K^a{}_b$ on pairs containing $\partial_t$.
	Together with transport of $F$ along the geodesics of the
	auxiliary connection, this gives
	\begin{equation}\label{eq:collar-evolution-system}
		\begin{aligned}
		\ddot\omega^\alpha+\omega^\alpha
		&=
		-\varepsilon^2
		\bigl((M_{11}+M_{22})\omega^\alpha
		-M_{\alpha\beta}\omega^\beta\bigr),\\
		\dot\Gamma^1{}_2
		&=
		\varepsilon^2(M_{02}\omega^1-M_{01}\omega^2),\\
		v^ae_a(F)
		-\Gamma^b{}_{ca}v^av^c\frac{\partial F}{\partial v^b}
		&=0.
		\end{aligned}
	\end{equation}
	In the last equation, $e_a(F)$ differentiates the base point while
	holding the velocity components fixed. The second term accounts
	for the change of those components along an auxiliary geodesic.
	Skew-symmetry of $\Gamma$ makes this vertical vector field tangent
	to $\Sph^2_v$, so the equation is well defined on unit directions.

	We prescribe the initial data
	\begin{equation}\label{eq:collar-initial-data}
		\begin{aligned}
		\omega^\alpha(0)&=\omega_\sigma^\alpha,
		&
		\dot\omega^\alpha(0)
		&=-\varepsilon A^\alpha{}_\beta\omega_\sigma^\beta,\\
		\Gamma^1{}_2(0)&=\Gamma^1{}_{2,\sigma},
		&
		F(0,q,v)&=\rho(q)\Phi(v),
		\end{aligned}
	\end{equation}
	where $\Gamma^1{}_{2,\sigma}$ is the round connection form and
	$\Phi$ is supplied by Lemma~\ref{lem:geodesic-representation}.

	The first equation in \eqref{eq:collar-evolution-system} is a
	perturbed Jacobi equation. More explicitly, set
	\begin{equation*}
		\widehat M
		=
		\begin{pmatrix}
		M_{22}&-M_{12}\\
		-M_{12}&M_{11}
		\end{pmatrix}.
	\end{equation*}
	Then
	\begin{equation*}
		\ddot C_\varepsilon+C_\varepsilon
		=
		-\varepsilon^2\widehat M C_\varepsilon,
		\qquad
		C_\varepsilon(0)=I,
		\qquad
		\dot C_\varepsilon(0)=-\varepsilon A.
	\end{equation*}
	Consequently,
	\begin{equation}\label{eq:collar-coframe-integral}
		C_\varepsilon(t)
		=
		\cos t\,I-\varepsilon\sin t\,A
		-\varepsilon^2
		\int_0^t
		\sin(t-s)\widehat M(s)C_\varepsilon(s)\,ds.
	\end{equation}
	This exhibits the desired metric as a correction of order
	$\varepsilon^2$ to the cosine collar. It is an integral equation:
	the matrix $\widehat M$ still depends on the unknown density $F$.

	\medskip
	\noindent
	\emph{Step 2: existence on a fixed collar.}\par\nopagebreak[4]

	We give the existence argument before checking the remaining
	geometric constraints. Its useful feature is that the auxiliary
	connection is an independent unknown. Thus the coefficients in
	the transport equation involve the coframe and the connection,
	but no spatial derivatives of these unknowns.

	Choose a compact interval
	\begin{equation*}
		J=[t_-,t_+]\subset(-\tau,t_2),
		\qquad
		t_-<0<t_1<t_+.
	\end{equation*}
	The round density $F_{\mathrm{rd}}$ from
	Lemma~\ref{lem:geodesic-representation} is defined on a
	neighborhood of this interval, and its support satisfies
	$v^0\geq\kappa$ for some $\kappa>0$.

	For the analytic argument, regard the velocity variable as
	belonging to the fixed sphere bundle
	$
	P=S(\mathbb R\oplus T\Sph^2),
	$
	with the norm induced by $1\oplus\sigma$. This is a compact
	four-dimensional manifold. At $(t,q)$, we identify its fiber
	with the $G_\varepsilon$-unit tangent sphere by
	\begin{equation*}
		(v^0,v^\top)\longmapsto
		v^0\partial_t+C_\varepsilon^{-1}v^\top,
		\qquad v^\top\in T_q\Sph^2,
	\end{equation*}
	where the second term is viewed as tangent to the slice
	$\{t\}\times\Sph^2$. We use the round case
	$C_0(t)=\cos t\,I$ to regard $F_{\mathrm{rd}}$ as a function
	on $J\times P$. Under time-independent changes of
	the local oriented orthonormal coframe, $C_\varepsilon$ and
	$\dot C_\varepsilon$ transform as bundle endomorphisms, while
	$
	\zeta=\Gamma^1{}_2-\Gamma^1{}_{2,\sigma}
	$
	is a globally defined one-form. We write
	$
	Z=(C_\varepsilon,\dot C_\varepsilon,\zeta).
	$

	On the transverse region, division of the transport equation by
	$v^0$ puts it in evolution form. To solve an equation on all of
	$P$, choose a smooth function $\vartheta(v^0)$ such that
	\begin{equation*}
		\vartheta(v^0)=\frac1{v^0}
		\quad\text{for }v^0\geq\kappa/2,
		\qquad
		\vartheta(v^0)=0
		\quad\text{for }v^0\leq\kappa/4.
	\end{equation*}
	Initially replace division by $v^0$ by multiplication by
	$\vartheta(v^0)$. The resulting equation is
	\begin{equation*}
		\partial_tF+B(Z)F=0,
	\end{equation*}
	where $B(Z)$ is a smooth vector field on $P$.
	As long as $C_\varepsilon$ is invertible, its coefficients depend
	smoothly on $Z$, $C_\varepsilon^{-1}$, and the fixed background
	data, without involving spatial derivatives of $Z$.

	For an integer $m\geq5$, use $H^m(\Sph^2)$ for $Z$ and
	$H^m(P)$ for $F$. Fiber integration gives a bounded map
	\begin{equation*}
		F\longmapsto M,
		\qquad
		H^m(P)\longrightarrow H^m(\Sph^2).
	\end{equation*}
	The transport energy estimate has the form
	\begin{equation*}
		\frac{d}{dt}\|F\|_{H^m}
		\leq
		C\left(
		\|B(Z)\|_{W^{1,\infty}}\|F\|_{H^m}
		+
		\|B(Z)\|_{H^m}\|F\|_{W^{1,\infty}}
		\right).
	\end{equation*}
	The Sobolev embedding theorem controls the first derivatives in
	$L^\infty$. Since the coefficients contain no spatial derivatives
	of $Z$, this estimate closes together with the evolution
	estimates for $Z$. Picard iteration, with uniform bounds in
	$H^m$ and convergence in $H^{m-1}$, therefore gives a unique
	local solution.

	To extend this solution over the fixed interval $J$, compare it
	with the round solution
	\begin{equation*}
		C_0(t)=\cos t\,I,
		\qquad
		\dot C_0(t)=-\sin t\,I,
		\qquad
		\zeta_0(t)=0,
		\qquad
		F_0=F_{\mathrm{rd}}.
	\end{equation*}
	The round coframe is uniformly invertible on $J$. Moreover, this
	reference solution satisfies the modified transport equation,
	because the modification is inactive on its support.

	Set
	\begin{equation*}
		E_m(t)
		=
		\|Z_\varepsilon(t)-Z_0(t)\|_{H^m(\Sph^2)}
		+
		\|F_\varepsilon(t)-F_{\mathrm{rd}}(t)\|_{H^m(P)}.
	\end{equation*}
	The initial data give $E_m(0)\leq C_m|\varepsilon|$.
	For the transport difference, we use
	\begin{equation*}
		\begin{aligned}
		\partial_t(F_\varepsilon-F_{\mathrm{rd}})
		+B(Z_\varepsilon)(F_\varepsilon-F_{\mathrm{rd}})
		=
		-\bigl(B(Z_\varepsilon)-B(Z_0)\bigr)
		F_{\mathrm{rd}}.
		\end{aligned}
	\end{equation*}
	The derivative on the right-hand side falls on the known smooth
	function $F_{\mathrm{rd}}$. On a fixed small neighborhood of the
	round solution in the $H^5$ norms, the tame energy estimates give
	\begin{equation*}
		\frac{d}{dt}E_m
		\leq C_m(E_m+\varepsilon^2).
	\end{equation*}
	The same estimate holds after reversing time. Gronwall's
	inequality and a continuation argument yield
	\begin{equation*}
		\sup_{t\in J}E_m(t)\leq C_m|\varepsilon|.
	\end{equation*}
	First apply this estimate at $m=5$ to choose a single sufficiently
	small parameter interval. It preserves the invertibility of
	$C_\varepsilon$ and gives a solution throughout $J$. Higher-order
	tame estimates are linear in the higher norms, with constants
	depending on $m$, the fixed smooth background, and the controlled
	low-order norms. They therefore propagate every Sobolev norm on
	the same interval, without further shrinking the parameter interval.

	The characteristic vector fields are uniformly close to the
	round characteristic field. Hence the characteristics issuing
	from the support of $\rho\Phi$ remain close to their round
	counterparts on $J$. After reducing $\varepsilon_0$, their
	velocity components satisfy
	$
	v^0\geq\kappa/2.
	$
	Thus the modification is inactive on the support of
	$F_\varepsilon$, and the solution satisfies the original
	transport equation in \eqref{eq:collar-evolution-system}.
	Constancy along characteristics also gives $F_\varepsilon\geq0$.

	Differentiating the equations with respect to $\varepsilon$,
	using higher spatial regularity at each differentiation, gives
	smooth dependence on $\varepsilon$. Finally, covariance under
	time-independent frame changes and uniqueness identify the local
	solutions on overlaps. We therefore obtain globally defined
	objects on the fixed open collar
$U=(t_-,t_+)\times\Sph^2.$

	\medskip
	\noindent
	\emph{Step 3: the conservation law for the moment tensor.}\par\nopagebreak[4]

	The transport equation implies a conservation law for $M$ that
	will preserve the curvature constraints. This calculation uses
	only the skew-symmetry of the auxiliary connection.

	At a fixed base point, the vertical transport vector field on
	the velocity sphere has components
	\begin{equation*}
		B^b(v)=-\Gamma^b{}_{ca}v^av^c.
	\end{equation*}
	It is tangent to $\Sph^2_v$ and satisfies
	\begin{equation*}
		\operatorname{div}_{\Sph^2_v}B
		=
		-\Gamma^a{}_{ca}v^c.
	\end{equation*}
	Multiplying the transport equation by $v^b$ and integrating over
	$\Sph^2_v$, we obtain
	\begin{equation}\label{eq:collar-moment-conservation}
		e_aM^{ab}
		+\Gamma^a{}_{ca}M^{cb}
		+\Gamma^b{}_{ca}M^{ac}
		=0.
	\end{equation}
	This is precisely the covariant divergence equation
	\begin{equation*}
		M^{ab}{}_{;a}=0
	\end{equation*}
	for the auxiliary connection. In particular, it holds before
	that connection has been shown to be Levi--Civita.

	\medskip
	\noindent
	\emph{Step 4: verification of the geometric constraints.}\par\nopagebreak[4]

	Define the torsion and curvature errors by
	\begin{equation*}
		\Theta^a
		=
		d\omega^a+\Gamma^a{}_b\wedge\omega^b,
		\qquad
		E^a{}_b
		=
		d\Gamma^a{}_b+\Gamma^a{}_c\wedge\Gamma^c{}_b
		-K^a{}_b.
	\end{equation*}
	The definition $\Gamma^\alpha{}_0=\dot\omega^\alpha$ implies
	$
	\iota_{\partial_t}\Theta^a=0.
	$
	Likewise, the first two evolution equations in
	\eqref{eq:collar-evolution-system} imply
	$
	\iota_{\partial_t}E^a{}_b=0.
	$
	Thus both errors are spatial two-forms, which we view as
	time-dependent two-forms on the fixed sphere $\Sph^2$.

	Symmetry of $M$ implies the algebraic Bianchi identity
	\begin{equation*}
		K^a{}_b\wedge\omega^b=0.
	\end{equation*}
	The Cartan identities therefore give
	\begin{equation}\label{eq:collar-error-identities}
		d^\Gamma\Theta^a=E^a{}_b\wedge\omega^b,
		\qquad
		d^\Gamma E=-d^\Gamma K.
	\end{equation}
	We explain how the conservation law removes the derivatives of
	$M$ from the second identity. Introduce the cyclic vector of
	two-forms
	\begin{equation*}
		\boldsymbol\beta
		=
		\begin{pmatrix}
		\omega^1\wedge\omega^2\\
		\omega^2\wedge\omega^0\\
		\omega^0\wedge\omega^1
		\end{pmatrix}.
	\end{equation*}
	Because the torsion forms are spatial,
	\begin{equation*}
		d^\Gamma\boldsymbol\beta
		=
		dt\wedge
		\begin{pmatrix}
		0\\
		\Theta^2\\
		-\Theta^1
		\end{pmatrix}.
	\end{equation*}
	Write the prescribed curvature as the cyclic vector
	\begin{equation*}
		\boldsymbol K=
		\begin{pmatrix}
		K^1{}_2\\
		K^2{}_0\\
		K^0{}_1
		\end{pmatrix}
		=
		(I+\varepsilon^2M)\boldsymbol\beta,
	\end{equation*}
	so the covariant product rule gives
	\begin{equation*}
		d^\Gamma\boldsymbol K
		=\varepsilon^2(d^\Gamma M)\wedge\boldsymbol\beta
		+(I+\varepsilon^2M)d^\Gamma\boldsymbol\beta.
	\end{equation*}
	By symmetry of $M$ and \eqref{eq:collar-moment-conservation},
	\begin{equation*}
		\bigl[(d^\Gamma M)\wedge\boldsymbol\beta\bigr]_a
		=M^{ab}{}_{;b}\,\omega^0\wedge\omega^1\wedge\omega^2
		=0.
	\end{equation*}
	Thus only the terms linear in the torsion $\Theta$ remain.

	Taking the components containing $dt$ in
	\eqref{eq:collar-error-identities} now gives
	\begin{equation}\label{eq:collar-error-evolution}
		\begin{aligned}
		\dot\Theta^a
		&=E^a{}_0,\\
		\dot E^\alpha{}_0
		&=
		-\Theta^\alpha
		-\varepsilon^2
		\bigl((M_{11}+M_{22})\Theta^\alpha
		-M_{\alpha\beta}\Theta^\beta\bigr),\\
		\dot E^1{}_2
		&=
		\varepsilon^2(M_{02}\Theta^1-M_{01}\Theta^2).
		\end{aligned}
	\end{equation}
	For each fixed $q$, these form a homogeneous linear system of
	ordinary differential equations for the error forms.

	We verify that their initial values vanish. The round surface
	connection and the symmetry of $A$ give
	$
	\Theta^a(0)=0.
	$
	The normalization of $\Phi$ gives
	\begin{equation*}
		M_{00}(0)=\rho,
		\qquad
		M_{0\alpha}(0)=0.
	\end{equation*}
	The Codazzi equation for $A$ then yields
	\begin{equation*}
		E^\alpha{}_0(0)=0.
	\end{equation*}
	For the tangential curvature component, the initial connection
	has
	\begin{equation*}
		\left.
		\bigl(d\Gamma^1{}_2+\Gamma^1{}_0\wedge\Gamma^0{}_2\bigr)
		\right|_{t=0}
		=
		\bigl(1-\varepsilon^2\det_\sigma A\bigr)\,d\mu_\sigma,
	\end{equation*}
	whereas the prescribed curvature form restricts to
	\begin{equation*}
		K^1{}_2\big|_{t=0}
		=
		(1+\varepsilon^2\rho)\,d\mu_\sigma.
	\end{equation*}
	Consequently,
	\begin{equation*}
		E^1{}_2(0)
		=
		-\varepsilon^2(\det_\sigma A+\rho)\,d\mu_\sigma
		=
		0.
	\end{equation*}

	Uniqueness for \eqref{eq:collar-error-evolution} therefore gives
	\begin{equation*}
		\Theta=0,
		\qquad
		E=0
		\qquad\text{throughout }U.
	\end{equation*}
	The auxiliary connection is thus metric-compatible and
	torsion-free, hence is the Levi--Civita connection of
	$G_\varepsilon$. Its curvature forms are exactly
	\eqref{eq:collar-prescribed-curvature}. This proves the exact
	identity
	\begin{equation*}
		T(G_\varepsilon)
		=
		\varepsilon^2M_{ab}\,\omega^a\otimes\omega^b
		\geq0,
	\end{equation*}
	and therefore $\sec_{G_\varepsilon}\geq1$ on $U$.

	\medskip
	\noindent
	\emph{Step 5: boundary data and the asymptotic expansions.}\par\nopagebreak[4]

	The initial coframe data immediately give
	\begin{equation*}
		\gamma_\varepsilon(0)=\sigma,
		\qquad
		\dot\gamma_\varepsilon(0)=-2\varepsilon A.
	\end{equation*}
	Since the outward unit normal at $t=0$ is $-\partial_t$, the
	boundary second fundamental form is
	\begin{equation*}
		h_{G_\varepsilon}
		=
		-\frac12\dot\gamma_\varepsilon(0)
		=
		\varepsilon A.
	\end{equation*}

	At $\varepsilon=0$, the metric and connection are round, and
	uniqueness for the transport equation gives
	$F_0=F_{\mathrm{rd}}$. Hence
	Lemma~\ref{lem:geodesic-representation} implies
	\begin{equation*}
		\left.
		M_{ab}\,\omega^a\otimes\omega^b
		\right|_{\varepsilon=0}
		=
		D.
	\end{equation*}
	Smooth dependence on $\varepsilon$ yields
	\begin{equation*}
		T(G_\varepsilon)
		=
		\varepsilon^2D+O(|\varepsilon|^3).
	\end{equation*}
	On the slab $[t_0,t_1]\times\Sph^2$, we have $D\geq c_0g_0$.
	Thus, for sufficiently small $|\varepsilon|$,
	\begin{equation*}
		T(G_\varepsilon)
		\geq
		\frac{c_0}{2}\varepsilon^2g_0
		\geq
		c\varepsilon^2G_\varepsilon
	\end{equation*}
	for some $c>0$. This is equivalent to the asserted strict
	sectional curvature bound.

	Finally, differentiating \eqref{eq:collar-coframe-integral} at
	$\varepsilon=0$ gives
	\begin{equation*}
		\left.\partial_\varepsilon C_\varepsilon\right|_{\varepsilon=0}
		=
		-\sin t\,A.
	\end{equation*}
	Therefore the first metric variation has no normal or mixed
	components, and its tangential part is
	\begin{equation*}
		\left.\partial_\varepsilon G_\varepsilon\right|_{\varepsilon=0}
		=
		-2\sin t\cos t\,A.
	\end{equation*}
	The vector field $V=f\partial_t-\tan t\,\nabla_\sigma f$ has the
	following expression along $X(t,q)=(\cos t\,q,\sin t)$ in
	$\mathbb R^3\times\mathbb R$:
	\begin{equation*}
		V(t,q)=\bigl(-\sin t\,F_f(q),\cos t\,f(q)\bigr).
	\end{equation*}
	For $Y\in T_q\Sph^2$, Euclidean differentiation and $dF_f=A$ give
	\begin{equation*}
		\partial_tV=\bigl(-\cos t\,F_f,-\sin t\,f\bigr),
		\qquad
		d_qV(Y)=\bigl(-\sin t\,A(Y),\cos t\,Yf\bigr).
	\end{equation*}
	Pairing these derivatives with
	$\partial_tX=(-\sin t\,q,\cos t)$ and
	$d_qX(Y)=(\cos t\,Y,0)$ yields
	\begin{equation*}
		\begin{aligned}
		(\operatorname{Lie}_Vg_0)(\partial_t,\partial_t)&=0,
		\qquad
		(\operatorname{Lie}_Vg_0)(\partial_t,Y)=0,
		\\
		(\operatorname{Lie}_Vg_0)|_{T\Sph^2}
		&=-2\sin t\cos t\,A.
		\end{aligned}
	\end{equation*}
	This proves the second expansion in
	\eqref{eq:collar-expansions}.
\end{proof}

\subsection{The curvature correction with prescribed support}

To extend the collar metric, we need to prescribe its quadratic
curvature term without changing the metric near the boundary.
The compatibility condition is the linearized Bianchi identity;
gnomonic coordinates will turn it into an ordinary Euclidean
divergence equation. Let
\begin{equation*}
	Lk
	=
	\left.\frac{d}{d\lambda}\right|_{\lambda=0}
	T(g_0+\lambda k)
\end{equation*}
be the linearization of $T$ at the round metric.
If a smooth extension $\bar g_\varepsilon$ satisfies
\begin{equation*}
	T(\bar g_\varepsilon)
	=
	\varepsilon^2E+O(|\varepsilon|^3),
\end{equation*}
then adding a metric correction $\varepsilon^2k$ gives
\begin{equation*}
	T(\bar g_\varepsilon+\varepsilon^2k)
	=
	\varepsilon^2(E+Lk)+O(|\varepsilon|^3).
\end{equation*}
Thus, to replace the quadratic curvature term $E$ by $D$, we
will solve $Lk=D-E$. The correction must be supported in the
interior cap so that the exact collar metric is preserved.

Naturality of $T$ and the identity $T(g_0)=0$
imply
\begin{equation}\label{eq:linearized-curvature-identities}
	L(\operatorname{Lie}_Wg_0)=0,
	\qquad
	\operatorname{div}_{g_0}(Lk)=0.
\end{equation}
The second identity follows by differentiating the contracted
Bianchi identity at $g_0$. It gives the necessary compatibility
condition for the linearized equation. All divergences in this
linearized problem are taken with respect to the fixed round
metric $g_0$.

\begin{lemma}\label{lem:supported-linear-correction}
	Let
	$
	\Omega=\{t>a\},
	$ and $0<a<\frac{\pi}{2}.$
	Suppose that $S$ is a smooth symmetric two-tensor on $H$ satisfying
	\begin{equation*}
		\operatorname{div}_{g_0}S=0,
		\qquad
		\operatorname{supp}S\subset\overline\Omega.
	\end{equation*}
	Then there exists a smooth symmetric two-tensor $k$ on $H$ such
	that
	\begin{equation*}
		Lk=S,
		\qquad
		\operatorname{supp}k\subset\overline\Omega.
	\end{equation*}
	In particular, $k$ is smooth across $\partial\Omega$ and vanishes
	to all orders there. The tensor $k$ can be chosen linearly in $S$.
\end{lemma}

The proof uses the gnomonic coordinates to turn the linearized
curvature equation into a Euclidean double-curl equation.
We construct its potential by integrating along outward rays
from a point inside the ball. This direction of integration
preserves support: a ray starting outside the ball and moving
away from an interior point never encounters the prescribed
tensor. The potential may be singular at the chosen center,
but averaging over the centers removes this singularity while
preserving both the equation and the support.

\begin{proof}
	In gnomonic coordinates, $\Omega$ corresponds to the Euclidean
	ball $B_R(0)$, where $R=\cot a$. Write
	$
	s=1+|y|^2.
	$
	For a Euclidean matrix field $p$, we use row-wise divergence
	and curl:
	\begin{equation*}
		(\operatorname{Div}p)_i=\partial_jp_{ij},
		\qquad
		(\operatorname{Curl}p)_{ij}
		=
		\epsilon_{jk\ell}\partial_kp_{i\ell},
	\end{equation*}
	where $\epsilon_{ijk}$ is the Euclidean alternating symbol.
	For brevity, set
	\begin{equation*}
		\operatorname{Inc}p
		=
		\operatorname{Curl}\bigl((\operatorname{Curl}p)^T\bigr).
	\end{equation*}
	Euclidean indices are raised and lowered using the standard
	Euclidean metric.

	\medskip
	\noindent
	\emph{Step 1: reduction to the Euclidean double-curl equation.}\par\nopagebreak[4]

	We claim that, for every smooth symmetric matrix field $p$,
	\begin{equation}\label{eq:linearized-curvature-double-curl}
		\bigl[L(s^{-1}p)\bigr]^{ij}
		=
		-\frac12s^3(\operatorname{Inc}p)^{ij}.
	\end{equation}
	Here $s^{-1}p$ denotes the covariant tensor with components
	$s^{-1}p_{ij}$, while the indices on the left-hand side are
	raised using $g_0$.

	We first compare the second-order terms. The linearization
	formulas for scalar and Ricci curvature give
	\begin{equation*}
		(Lh)^{ij}
		=
		-\frac{1}{2\det g_0}
		\epsilon_{imn}\epsilon_{jk\ell}
		\partial_m\partial_k h_{n\ell}
		+
		\text{terms involving at most first derivatives of }h.
	\end{equation*}
	Since $\det g_0=s^{-4}$, the second-order part of $L(s^{-1}p)$
	is $-\frac12s^3\operatorname{Inc}p$. Hence
	\begin{equation*}
		R_0(p)^{ij}
		=
		\bigl[L(s^{-1}p)\bigr]^{ij}
		+\frac12s^3(\operatorname{Inc}p)^{ij}
	\end{equation*}
	is a differential operator of order at most one.

	We show that this remaining operator vanishes by using
	diffeomorphism invariance. The gnomonic Christoffel symbols give,
	for every one-form $\alpha$,
	\begin{equation*}
		\nabla_i\alpha_j+\nabla_j\alpha_i
		=
		s^{-1}
		\bigl[\partial_i(s\alpha_j)+\partial_j(s\alpha_i)\bigr].
	\end{equation*}
	Thus, if $\beta$ is a Euclidean one-form and
	$\alpha=s^{-1}\beta$, then
	\begin{equation*}
		\operatorname{Lie}_Wg_0
		=
		s^{-1}\operatorname{sym}\partial\beta,
		\qquad
		W=\alpha^{\sharp_{g_0}},
	\end{equation*}
	where
	\begin{equation*}
		(\operatorname{sym}\partial\beta)_{ij}
		=
		\partial_i\beta_j+\partial_j\beta_i.
	\end{equation*}
	By \eqref{eq:linearized-curvature-identities},
	\begin{equation*}
		L\bigl(s^{-1}\operatorname{sym}\partial\beta\bigr)=0.
	\end{equation*}
	Commutativity of Euclidean partial derivatives also gives
	\begin{equation*}
		\operatorname{Inc}(\operatorname{sym}\partial\beta)=0.
	\end{equation*}
	Therefore $R_0$ annihilates every Euclidean symmetric
	gradient.

	The value and first derivatives of a symmetric gradient can
	be prescribed arbitrarily at a point. Indeed, for any prescribed
	symmetric first jet $(p_{ij},p_{ij,k})$, choose
	\begin{equation*}
		\partial_j\beta_i=\frac12p_{ij},
		\qquad
		\partial_j\partial_k\beta_i
		=
		\frac12\bigl(p_{ij,k}+p_{ik,j}-p_{jk,i}\bigr)
	\end{equation*}
	at that point. These choices give
	\begin{equation*}
		(\operatorname{sym}\partial\beta)_{ij}=p_{ij},
		\qquad
		\partial_k(\operatorname{sym}\partial\beta)_{ij}=p_{ij,k}.
	\end{equation*}
	Since $R_0$ has order at most one and annihilates all
	such jets, it is identically zero. This proves
	\eqref{eq:linearized-curvature-double-curl}, including its
	lower-order terms.

	Now define the Euclidean symmetric tensor
	\begin{equation}\label{eq:weighted-linearized-source}
		U^{ij}=s^{-3}S^{ij},
	\end{equation}
	where the indices on $S$ are raised using $g_0$.
	The divergence identity in Lemma~\ref{lem:divergence-free-tensor}
	gives
	$
	\partial_iU^{ij}=0.
	$
	Since $U$ is symmetric, this is equivalent to
	$\operatorname{Div}U=0$. Moreover, $U$ is smooth on $\mathbb R^3$
	and supported in $\overline{B_R(0)}$.

	By \eqref{eq:linearized-curvature-double-curl}, it remains to
	construct a smooth symmetric matrix field $p$, supported in
	$\overline{B_R(0)}$, such that
	$
		\operatorname{Inc}p=U.
	$
	We will then take $k=-2s^{-1}p$.

	\medskip
	\noindent
	\emph{Step 2: a vector potential obtained from outward rays.}\par\nopagebreak[4]

	Fix $b\in B_R(0)$. For a smooth, compactly supported,
	divergence-free Euclidean vector field $v$, define, for $y\neq b$,
	\begin{equation}\label{eq:outward-ray-vector-potential}
		K_bv(y)
		=
		\int_1^\infty
		u\,(y-b)\times v\bigl(b+u(y-b)\bigr)\,du.
	\end{equation}
	This integral samples $v$ along the ray starting at $y$ and
	continuing away from $b$.

	Put $x=y-b$. A direct calculation using $\operatorname{div}v=0$
	gives
	\begin{equation*}
		\operatorname{curl}_y
		\bigl[u\,x\times v(b+ux)\bigr]
		=
		-\frac{d}{du}\bigl[u^2v(b+ux)\bigr].
	\end{equation*}
	Since $v$ is compactly supported, integration in $u$ yields
	\begin{equation*}
		\operatorname{curl}K_bv=v
		\qquad\text{on }\mathbb R^3\setminus\{b\}.
	\end{equation*}

	This construction also preserves support in the ball. Indeed,
	if $v$ is supported in $\overline{B_R(0)}$ and
	$y\notin\overline{B_R(0)}$, then
	\begin{equation*}
		b+u(y-b)\notin\overline{B_R(0)}
		\qquad\text{for every }u\geq1,
	\end{equation*}
	by convexity of the ball. Hence $K_bv(y)=0$.

	\medskip
	\noindent
	\emph{Step 3: a symmetric double-curl potential.}\par\nopagebreak[4]

	Apply $K_b$ to each row of $U$, and write
	$
	Q_b=K_bU.
	$
	Then
	\begin{equation*}
		\operatorname{Curl}Q_b=U
		\qquad\text{away from }b.
	\end{equation*}
	For $x=y-b$, let $J_x$ denote the skew-symmetric matrix
	$
	J_xw=x\times w.$ In matrix notation,
	\begin{equation*}
		Q_b(y)
		=
		\int_1^\infty
		u\,U(b+ux)J_x^T\,du.
	\end{equation*}
	Since $U$ is symmetric and $J_x$ is skew-symmetric,
	$
	\operatorname{tr}Q_b=0.
	$

	The symmetry of $\operatorname{Curl}Q_b=U$ gives another useful
	identity:
	\begin{equation*}
		\partial_j(Q_b)_{ji}
		=
		\partial_i\operatorname{tr}Q_b
		=
		0.
	\end{equation*}
	Thus every row of $Q_b^T$ is divergence-free, and we may apply
	the vector potential construction a second time. Set
	\begin{equation*}
		P_bU=K_b(Q_b^T).
	\end{equation*}
	The tensors involved are smooth away from $b$ and vanish outside
	the closed ball, so the outward-ray integrals are well defined
	there. We obtain
	\begin{equation*}
		\operatorname{Curl}(P_bU)=Q_b^T,
		\qquad
		\operatorname{Inc}(P_bU)=U
		\qquad\text{away from }b.
	\end{equation*}

	Combining the two radial integrals gives the explicit formula
	\begin{equation}\label{eq:outward-ray-double-curl-potential}
		P_bU(y)
		=
		\int_1^\infty
		u(u-1)\,
		J_{y-b}\,
		U\bigl(b+u(y-b)\bigr)\,
		J_{y-b}^T\,du.
	\end{equation}
	Indeed, the two radial parameters $r,\lambda\geq1$ combine as
	$u=r\lambda$; the factor
	$r^2\lambda\,d\lambda\,dr$ becomes $u\,du\,dr$, with
	$1\leq r\leq u$. Integration in $r$ gives $u(u-1)$.
	Formula \eqref{eq:outward-ray-double-curl-potential} shows directly
	that $P_bU$ is symmetric.

	We must also check the equations at the center $b$.
	For $y$ near $b$, compact support of $U$ restricts the radial
	parameter to an interval of length $O(|y-b|^{-1})$.
	The preceding formulas give
	\begin{equation*}
		Q_b(y)=O(|y-b|^{-1}),
		\qquad
		P_bU(y)=O(|y-b|^{-1}).
	\end{equation*}
	These singularities are locally integrable. When either
	first-curl identity is tested by integration by parts outside
	$B_\delta(b)$, the boundary term is bounded by
	\begin{equation*}
		C\delta^2\delta^{-1}=C\delta,
	\end{equation*}
	and therefore tends to zero as $\delta\downarrow0$.
	Consequently,
	\begin{equation*}
		\operatorname{Curl}Q_b=U,
		\qquad
		\operatorname{Curl}(P_bU)=Q_b^T
	\end{equation*}
	hold distributionally on all of $\mathbb R^3$. In particular,
	\begin{equation*}
		\operatorname{Inc}(P_bU)=U
	\end{equation*}
	holds distributionally across the center as well.

	\medskip
	\noindent
	\emph{Step 4: averaging the centers to obtain a smooth potential.}\par\nopagebreak[4]

	Choose a nonnegative function
	\begin{equation*}
		\theta\in C_c^\infty(B_R(0)),
		\qquad
		\int_{\mathbb R^3}\theta(b)\,db=1,
	\end{equation*}
	and define
	\begin{equation}\label{eq:averaged-double-curl-potential}
		PU(y)
		=
		\int_{\mathbb R^3}\theta(b)P_bU(y)\,db.
	\end{equation}
	The singularity $|y-b|^{-1}$ is integrable in $b$, so this
	average is well defined. It is symmetric, and the distributional
	identities above give
	\begin{equation*}
		\operatorname{Inc}(PU)=U.
	\end{equation*}
	Each $P_bU$ is supported in $\overline{B_R(0)}$, and therefore
	so is $PU$.

	We now verify that the averaging removes the singularity
	completely. In
	\eqref{eq:averaged-double-curl-potential} and
	\eqref{eq:outward-ray-double-curl-potential}, make the changes
	of variables
	\begin{equation*}
		w=(u-1)(y-b),
		\qquad
		r=(u-1)^{-1}.
	\end{equation*}
	They give
	\begin{equation}\label{eq:averaged-potential-smoothness}
		PU(y)
		=
		\int_{\mathbb R^3}
		J_wU(y+w)J_w^T
		\left[
		\int_0^\infty
		(r+r^2)\theta(y-rw)\,dr
		\right]dw.
	\end{equation}

	Fix a compact set of $y$-values. On the support of the integrand,
	both $w$ and $r|w|$ are bounded, because $U$ and $\theta$ have
	compact support. For $w\neq0$, this gives
	\begin{equation*}
		|w|^2
		\int_0^{C/|w|}(r+r^2)\,dr
		\leq C'\bigl(1+|w|^{-1}\bigr).
	\end{equation*}
	The right-hand side is locally integrable in $\mathbb R^3$.

	With $w$ fixed, every derivative in $y$ in
	\eqref{eq:averaged-potential-smoothness} falls only on $U$ or
	$\theta$. In particular, differentiating $\theta(y-rw)$ in $y$
	introduces no additional factor of $r$. Thus every
	$y$-derivative of the integrand satisfies the same type of
	integrable bound, with a constant depending on the derivative
	order. Differentiation under the integral proves
	\begin{equation*}
		PU\in C^\infty(\mathbb R^3).
	\end{equation*}

	We have therefore constructed a smooth symmetric tensor
	satisfying
	\begin{equation*}
		\operatorname{Inc}(PU)=U,
		\qquad
		\operatorname{supp}(PU)\subset\overline{B_R(0)}.
	\end{equation*}
	Since it is globally smooth and vanishes outside the closed
	ball, it vanishes to all orders on the boundary sphere.

	\medskip
	\noindent
	\emph{Step 5: return to the spherical metric.}\par\nopagebreak[4]

	Define the covariant tensor $k$ in gnomonic coordinates by
	\begin{equation*}
		k_{ij}=-2s^{-1}(PU)_{ij}.
	\end{equation*}
	Using \eqref{eq:linearized-curvature-double-curl} and
	\eqref{eq:weighted-linearized-source}, we obtain
	\begin{equation*}
		(Lk)^{ij}
		=
		s^3\bigl(\operatorname{Inc}(PU)\bigr)^{ij}
		=
		s^3U^{ij}
		=
		S^{ij}.
	\end{equation*}
	Thus $Lk=S$.

	The support of $k$ is contained in the closed cap
	$\overline\Omega$. Its smoothness across $\partial\Omega$
	follows from the global smoothness of $PU$. Extending it by
	zero toward the equator gives a smooth tensor on $H$.
	Finally, with $\theta$ fixed, every operation in the construction
	is linear in $S$, so the resulting choice of $k$ is linear.
\end{proof}

\subsection{Gluing and the global curvature bound}

The exact collar controls the boundary region, while the strict
positivity of $D$ absorbs the remaining error on the interior cap.
The supported correction joins these two parts without changing
the prescribed boundary geometry.

\begin{figure}[htbp]
	\centering
	\resizebox{\linewidth}{!}{%
	\begin{tikzpicture}[
		x=1cm,
		y=1cm,
		font=\small,
		region/.style={draw=black!55, line width=0.5pt},
		guide/.style={densely dashed, draw=black!30},
		annotation/.style={align=center, inner sep=3pt}
		]


		\foreach \x in {0,1.2,3.1,5.4,7.1,12}{
		\draw[guide] (\x,0) -- (\x,3.4);
		}

		\draw[region, fill=blue!10]
		(0,2.4) rectangle (3.1,3.4);
		\draw[region, fill=orange!18]
		(3.1,2.4) rectangle (5.4,3.4);
		\draw[region, fill=black!5]
		(5.4,2.4) rectangle (12,3.4);

		\node[anchor=west, fill=white, inner sep=1pt]
		at (0,3.65)
		{$\bar g_\varepsilon$: preliminary extension};

		\node[annotation] at (1.55,2.9)
		{$G_\varepsilon$};

		\node[annotation] at (4.25,2.9)
		{cutoff\\transition};

		\node[annotation] at (8.7,2.9)
		{$g_0+\varepsilon\operatorname{Lie}_Vg_0$\\[2pt]
		{\footnotesize equal to $g_0$ near the pole}};

		\draw[region, fill=blue!10]
		(0,0.8) rectangle (3.1,1.8);
		\draw[region, fill=green!12]
		(3.1,0.8) rectangle (12,1.8);

		\node[anchor=west, fill=white, inner sep=1pt]
		at (0,2.05)
		{$g_\varepsilon$: final metric};

		\node[annotation] at (1.55,1.3)
		{$g_\varepsilon=G_\varepsilon$};

		\node[annotation] at (7.55,1.3)
		{$g_\varepsilon=\bar g_\varepsilon+\varepsilon^2k$\\[2pt]
		{\footnotesize
		$\operatorname{supp}k\subset\{t\geq a\}$}};

		\draw[
		decorate,
		decoration={brace, amplitude=5pt},
		line width=0.6pt
		]
		(1.2,4.15) -- (12,4.15)
		node[midway, above=8pt, align=center]
		{$D\geq c_0g_0$ on the interior cap $\{t\geq t_0\}$};

		\draw[->, line width=0.7pt]
		(-0.2,0) -- (12.55,0)
		node[right] {$t$};

		\foreach \x in {0,1.2,3.1,5.4,7.1,12}{
		\draw[line width=0.6pt]
		(\x,0.07) -- (\x,-0.07);
		}

		\node[below=5pt, align=center] at (0,0)
		{$0$\\{\footnotesize equator}};
		\node[below=5pt] at (1.2,0) {$t_0$};
		\node[below=5pt] at (3.1,0) {$a$};
		\node[below=5pt] at (5.4,0) {$b$};
		\node[below=5pt] at (7.1,0) {$t_1$};
		\node[below=5pt, align=center] at (12,0)
		{$\pi/2$\\{\footnotesize north pole}};

		\draw[
		decorate,
		decoration={brace, mirror, amplitude=5pt},
		line width=0.6pt
		]
		(0,-1.05) -- (7.1,-1.05)
		node[midway, below=8pt, align=center]
		{Exact collar metric available on $0\leq t\leq t_1$};

	\end{tikzpicture}%
	}
	\caption{
	A schematic figure for the construction
	}
	\label{fig:collar-extension-correction}
\end{figure}
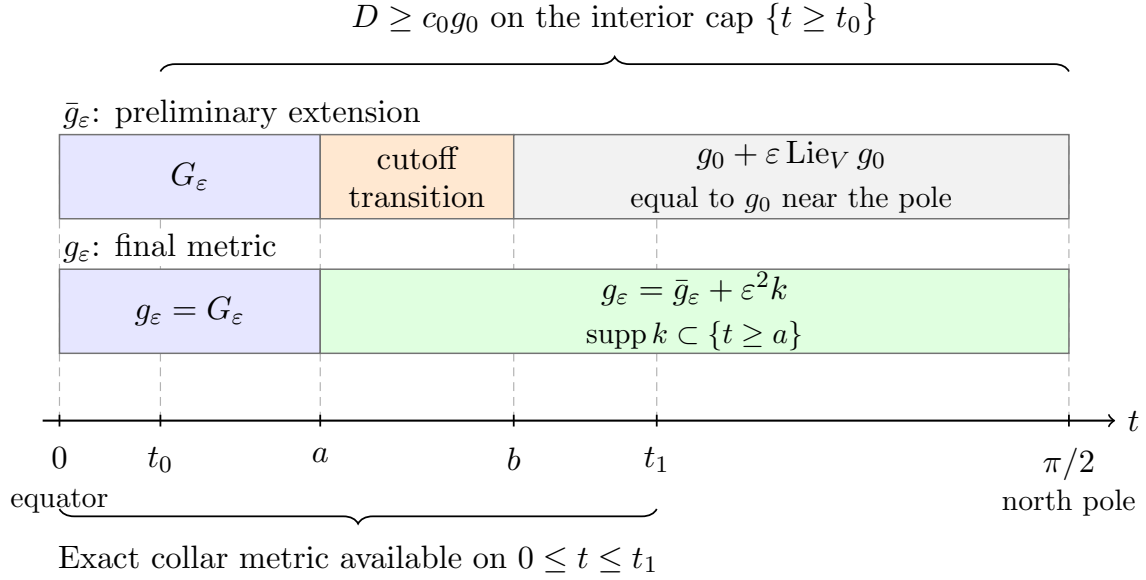

Figure~\ref{fig:collar-extension-correction} summarizes the regions used in the proof
of Theorem~\ref{thm:three-dimensional-counterexample}.

\begin{proof}[Proof of Theorem~\ref{thm:three-dimensional-counterexample}]
	Let $D$ be the tensor constructed in
	Lemma~\ref{lem:divergence-free-tensor}, and let
	$G_\varepsilon$ be the collar metrics supplied by
	Lemma~\ref{lem:exact-collar}. We extend these metrics to $H$
	and then correct their quadratic curvature term inside the
	cap, keeping the metric unchanged near the boundary.

	Recall that $t_0<a<b<t_1$. Choose a smooth cutoff function
	$\chi$ on $H$ such that
	\begin{equation*}
		0\leq\chi\leq1,
		\qquad
		\chi=1\ \text{on }\{t\leq a\},
		\qquad
		\chi=0\ \text{on }\{t\geq b\}.
	\end{equation*}
	Extend the vector field
	$
	V=f\partial_t-\tan t\,\nabla_\sigma f
	$
	smoothly to $H$, retaining its original values on a neighborhood
	of $\{t\leq b\}$ and making it vanish near the north pole.
	We denote the extension also by $V$.

	Define
	\begin{equation}\label{eq:preliminary-global-metric}
		\bar g_\varepsilon
		=
		g_0+\varepsilon\operatorname{Lie}_Vg_0
		+
		\chi\bigl(
		G_\varepsilon-g_0-\varepsilon\operatorname{Lie}_Vg_0
		\bigr),
	\end{equation}
	where the term multiplied by $\chi$ is extended by zero outside
	the collar. This defines a smooth family of symmetric tensors
	on $H$. It agrees with $G_\varepsilon$ on $\{t\leq a\}$ and
	with $g_0$ near the north pole.

	The first-variation formula in Lemma~\ref{lem:exact-collar}
	shows that
	\begin{equation*}
		G_\varepsilon-g_0-\varepsilon\operatorname{Lie}_Vg_0
		=
		O(\varepsilon^2)
	\end{equation*}
	on the support of $\chi$. Consequently,
	\begin{equation*}
		\bar g_\varepsilon
		=
		g_0+\varepsilon\operatorname{Lie}_Vg_0+O(\varepsilon^2)
	\end{equation*}
	on $H$. In particular, $\bar g_\varepsilon$ is positive definite
	for sufficiently small $|\varepsilon|$.

	Since $T(g_0)=0$ and
	$L(\operatorname{Lie}_Vg_0)=0$, the curvature expansion of
	$\bar g_\varepsilon$ starts at second order:
	\begin{equation}\label{eq:preliminary-curvature-expansion}
		T(\bar g_\varepsilon)
		=
		\varepsilon^2E+O(|\varepsilon|^3),
	\end{equation}
	where $E$ is a smooth symmetric tensor on $H$, independent of
	$\varepsilon$. The exact agreement with $G_\varepsilon$ on
	$\{t\leq a\}$ and the expansion
	\eqref{eq:collar-expansions} imply
	\begin{equation*}
		E=D
		\qquad\text{on }\{t\leq a\}.
	\end{equation*}

	We next verify the compatibility condition for the interior
	correction. The contracted Bianchi identity gives
	$
	\operatorname{div}_{\bar g_\varepsilon}
	T(\bar g_\varepsilon)=0.
	$
	Using \eqref{eq:preliminary-curvature-expansion} and
	$\bar g_\varepsilon=g_0+O(\varepsilon)$, we obtain
	\begin{equation*}
		0
		=
		\varepsilon^2\operatorname{div}_{g_0}E
		+
		O(|\varepsilon|^3).
	\end{equation*}
	Dividing by $\varepsilon^2$ and letting $\varepsilon\to0$ yields
	$
	\operatorname{div}_{g_0}E=0.
	$
	Therefore
	$
	S:=D-E
	$
	is a smooth symmetric tensor satisfying
	\begin{equation*}
		\operatorname{div}_{g_0}S=0,
		\qquad
		\operatorname{supp}S\subset\{t\geq a\}.
	\end{equation*}

	Apply Lemma~\ref{lem:supported-linear-correction} with
	$\Omega=\{t>a\}$. It gives a smooth symmetric tensor $k$ on $H$
	such that
	\begin{equation*}
		Lk=D-E,
		\qquad
		\operatorname{supp}k\subset\{t\geq a\}.
	\end{equation*}
	The tensor $k$ is fixed, independent of $\varepsilon$. Define
	the final metrics by
	$
	g_\varepsilon
	=
	\bar g_\varepsilon+\varepsilon^2k.
	$
	They form a smooth family and converge to $g_0$ as
	$\varepsilon\to0$, so they are positive definite for sufficiently
	small $|\varepsilon|$. Moreover,
	\begin{equation*}
		g_{\varepsilon=0}=g_0,
		\qquad
		g_\varepsilon=G_\varepsilon
		\quad\text{on }\{t\leq a\}.
	\end{equation*}

	Taylor expansion at the round metric gives
	\begin{equation*}
		T(g_\varepsilon)
		=
		T(\bar g_\varepsilon)
		+
		\varepsilon^2Lk
		+
		O(|\varepsilon|^3).
	\end{equation*}
	Indeed, the linearization at $\bar g_\varepsilon$ differs from
	$L$ by $O(|\varepsilon|)$, since
	$\bar g_\varepsilon=g_0+O(\varepsilon)$. Hence
	\begin{equation}\label{eq:final-curvature-expansion}
		\begin{aligned}
		T(g_\varepsilon)
		&=
		\varepsilon^2(E+Lk)+O(|\varepsilon|^3)\\
		&=
		\varepsilon^2D+O(|\varepsilon|^3).
		\end{aligned}
	\end{equation}

	We verify the curvature bound separately on the retained collar
	and the interior cap. On $\{t\leq a\}$, the metric is exactly
	$G_\varepsilon$, and Lemma~\ref{lem:exact-collar} gives
	$
	\sec_{g_\varepsilon}\geq1.
	$
	On $\{t\geq a\}$, the inequality $a>t_0$ and
	Lemma~\ref{lem:divergence-free-tensor} give
	$
	D\geq c_0g_0.
	$
	The remainder in \eqref{eq:final-curvature-expansion} is uniform
	on this compact cap. Thus, as quadratic forms,
	\begin{equation*}
		T(g_\varepsilon)
		\geq
		\varepsilon^2\bigl(c_0-C|\varepsilon|\bigr)g_0
		\geq0
	\end{equation*}
	for sufficiently small $|\varepsilon|$. By the
	three-dimensional characterization of the sectional curvature
	bound through $T$, this proves
	\begin{equation*}
		\sec_{g_\varepsilon}\geq1
		\qquad\text{on all of }H.
	\end{equation*}

	Finally, the metric is unchanged from $G_\varepsilon$ near the
	boundary. Therefore
	$
	g_\varepsilon|_{T\partial H}=\sigma.
	$
	The outward unit normal is $-\partial_t$, so
	\begin{equation*}
		h_{g_\varepsilon}
		=
		-\frac12\dot\gamma_\varepsilon(0)
		=
		\varepsilon A.
	\end{equation*}
	Since $A$ is Codazzi with respect to $\sigma$, the boundary
	second fundamental form is Codazzi for every $\varepsilon$.
	Moreover, $\rho\not\equiv0$ implies $A\not\equiv0$.
	Thus the boundary is not totally geodesic when
	$\varepsilon\neq0$, and $(H,g_\varepsilon)$ cannot be isometric
	to the standard hemisphere.
\end{proof}

\section{An inradius estimate in dimension three}
\label{sec:dimension-three}

The equality case of Theorem~\ref{thm:inradius-three} again reduces
to a collapsed great circle.
Distance comparison produces a nonnegative boundary excess whose
zero set is the contact set. We isolate the resulting spherical
argument in Proposition~\ref{prop:three-contact-rigidity}.
A contact barrier and the Codazzi equation turn a suitable regular
level curve of a scalar projection of the support map into a complete
great circle on which the map is constant. The G\'alvez--Mira theorem
supplies two distinct support values so that this level curve can be chosen;
Lemma~\ref{lem:circle-support-rigidity} then forces the tensor to vanish.
Recall that
\begin{equation*}
	r_{\mathrm{in}}(M)
	:=
	\max_{p\in M} d_g(p,\partial M).
\end{equation*}

Identify $\partial M$ with the unit round sphere $S^2$, and write
$\bar g=g|_{T\partial M}$. Denote its Levi-Civita connection by $\nabla$
and the ambient connection by $\nabla^M$. With the outward unit normal
$\nu$, set
\begin{equation*}
	h(X,Y)=g(\nabla^M_X\nu,Y),
	\qquad
	h(X,Y)=\bar g(AX,Y).
\end{equation*}
By hypothesis, $A$ is a self-adjoint Codazzi tensor.

The Gauss equation gives
\begin{equation}
	\label{eq:three-gauss}
	\det A
	=
	1-\sec_M(T\partial M)
	\leq 0.
\end{equation}

\subsection{Distance comparison at the boundary}
\label{sec:three-distance-comparison}

The following comparison holds in every dimension. Besides bounding
the second fundamental form at a nearest boundary point, it ensures
that the distance function is smooth there. This smoothness will
allow us to use the distance excess as a local barrier.

\begin{lemma}
	\label{lem:three-distance}
	Let $(M^n,g)$, $n\geq2$, be a compact connected Riemannian
	manifold with smooth nonempty boundary and $\sec_g\geq 1$.
	For $p\in\operatorname{int}M$, set
	\begin{equation*}
		R=d_g(p,\partial M),
		\qquad
		C_p=\{q\in\partial M:d_g(p,q)=R\}.
	\end{equation*}
	Then
	$R<\pi$.
	Moreover, the distance function $r=d_g(p,\cdot)$ is smooth
	in a neighborhood of $C_p$ and
	\begin{equation*}
		\nabla^M r(q)=\nu(q),
		\qquad
		h_q\leq \cot R\,\bar{g}_q
		\qquad\text{for every }q\in C_p.
	\end{equation*}
\end{lemma}

\begin{proof}
	Fix $q\in C_p$, and let
	$
	\alpha:[0,R]\longrightarrow M
	$
	be a minimizing unit-speed geodesic from $p$ to $q$.
	Its interior lies in $\operatorname{int}M$, and the first
	variation formula gives
	$
	\alpha'(R)=\nu(q).
	$

	For a vector field $V$ along $\alpha$ perpendicular to $\alpha'(s)$ satisfying
	\begin{equation*}
		V(0)=0,\qquad V(R)\in T_q\partial M,
	\end{equation*}
	the second variation formula gives
	\begin{equation*}
		Q(V,V)
		:=
		I_\alpha(V,V)-h_q(V(R),V(R))
		\geq 0,
	\end{equation*}
	where
	\begin{equation*}
		I_\alpha(V,W)
		=
		\int_0^R
		\left(
		\langle D_sV,D_sW\rangle
		-
		\langle R^M(V,\alpha') \alpha',W\rangle
		\right)\,ds.
	\end{equation*}

	We first show that $q$ is not conjugate to $p$ along $\alpha$.
	Otherwise, there would exist a nonzero perpendicular Jacobi
	field $J$ satisfying
	$J(0)=J(R)=0$ and $Q(J,J)=0$. Since $Q$ is positive semidefinite,
	$
	Q(J,V)=0
	$
	for every admissible $V$. Integration by parts gives
	\begin{equation*}
		0=Q(J,V)=\langle D_sJ(R),V(R)\rangle.
	\end{equation*}
	Since $D_sJ(R)$ is perpendicular to $\nu(q)$, it belongs to
	$T_q\partial M$. The arbitrariness of $V(R)$ therefore implies
	$D_sJ(R)=0$. Together with $J(R)=0$, this forces
	$J\equiv 0$, a contradiction.

	There are no conjugate points in $(0,R)$, since $\alpha$ is
	minimizing. Conjugate point comparison under $\sec_g\geq 1$
	therefore yields $R<\pi$.

	Every minimizing geodesic from $p$ to $q$ has terminal
	velocity $\nu(q)$. Uniqueness for the geodesic equation
	thus implies that $\alpha$ is the unique minimizing
	geodesic from $p$ to $q$.

	To justify smoothness of the distance at boundary points,
	take a smooth complete extension of $(M,g)$. Any curve
	from $p$ to $\partial M$ must first reach the boundary
	after a segment contained in $M$ of length at least $R$.
	Consequently, the distance from $p$ to $q$ in the extension
	is also $R$, and every minimizing geodesic to $q$ remains
	in $M$ until its endpoint. Thus $q$ is joined to $p$ by a
	unique minimizing geodesic and is not conjugate to $p$
	in the extension. The ambient distance function is
	therefore smooth near $q$.

	By transversality at $q$, the nearby minimizing geodesics
	to points of $M$ also remain in $M$. Hence the ambient
	distance agrees locally with $d_g(p,\cdot)$ on $M$.
	This proves smoothness near $C_p$, up to the boundary,
	and gives $\nabla^M r(q)=\nu(q)$.

	Finally, let $X\in T_q\partial M$, and let $P(s)$ be its
	parallel transport along $\alpha$. Taking
	$
	V(s)=\frac{\sin s}{\sin R}\,P(s),
	$
	the second variation inequality and $\sec_g\geq 1$ then imply
	\begin{equation*}
		\begin{aligned}
		h_q(X,X)
		&\leq I_\alpha(V,V)\\
		&\leq
		\frac{|X|^2}{\sin^2 R}
		\int_0^R(\cos^2s-\sin^2s)\,ds\\
		&=\cot R\,|X|^2.
		\end{aligned}
	\end{equation*}
	The claimed tensor inequality follows.
\end{proof}

\subsection{Contact barriers and kernel propagation}
\label{sec:three-contact-propagation}

We need two consequences of the distance inequality. The first
prevents a separated part of the contact set from avoiding a closed
geodesic semicircle. A subset of $C$ is called relatively clopen if
it is both open and closed in $C$.

\begin{lemma}
	\label{lem:three-convex-barrier}
	Let $U\subset S^2$ be open, and let $u\in C^\infty(U)$ be
	nonnegative. Suppose that
	\begin{equation*}
		C=\{x\in U:u(x)=0\}
	\end{equation*}
	is nonempty and compact, and that
	\begin{equation}
		\label{eq:three-no-positive-hessian}
		\lambda_{\min}(\nabla^2u)\leq 0
		\qquad\text{on }U.
	\end{equation}
	Then every nonempty relatively clopen subset of $C$
	intersects every closed geodesic semicircle of $S^2$.
\end{lemma}

\begin{proof}
	Let $K\subset C$ be nonempty and relatively clopen, and suppose
	that $K$ is disjoint from a closed geodesic semicircle $\Gamma$.

	By the convex supporting theorem of Jost--Xin--Yang
	\cite[Theorem~2.1]{JostXinYang}, the complement
	$S^2\setminus\Gamma$ admits a smooth strictly convex function
	on a neighborhood of any compact subset. Thus there exist an
	open neighborhood $W$ of $K$ and a smooth function $\psi$ on $W$
	such that
	\begin{equation*}
		\nabla^2\psi>0.
	\end{equation*}

	Since $K$ is relatively clopen in the compact set $C$, we can
	choose an open set $V$ satisfying
	\begin{equation*}
		K\subset V,\qquad
		\overline V\subset U\cap W,\qquad
		\overline V\cap C=K.
	\end{equation*}
	In particular, $u>0$ on $\partial V$. Set
	$
	m=\max_K\psi.
	$
	For sufficiently small $\varepsilon>0$, the function
	\begin{equation*}
		w=u-\varepsilon(\psi-m)
	\end{equation*}
	is positive on $\partial V$, while $w=0$ at a point of $K$
	where $\psi=m$. Therefore $w$ attains a nonpositive minimum
	at an interior point $x\in V$. At this point,
	\begin{equation*}
		\nabla^2u(x)
		\geq
		\varepsilon\nabla^2\psi(x)
		>0,
	\end{equation*}
	contradicting \eqref{eq:three-no-positive-hessian}.
\end{proof}

The second consequence is local propagation along the kernel at a
rank-one contact point. A kernel vector alone need not extend to a
curve in a support-map fiber when rank-two points occur nearby.
The barrier inequality below supplies this missing propagation.

\begin{lemma}
	\label{lem:three-rank-one-propagation}
	Let $h$ be a smooth symmetric $(0,2)$ Codazzi tensor on an open subset
	$U\subset S^2$, and let $A$ be its associated self-adjoint
	endomorphism. Assume that $\det A\leq 0$ on $U$.

	Suppose that $u,b\in C^\infty(U)$ satisfy
	\begin{equation*}
		u\geq 0,\qquad b>0,\qquad
		\nabla^2u+b h\leq 0.
	\end{equation*}
	If $u(q)=0$ and $A(q)\neq 0$, then there exists a unit-speed
	geodesic segment $\sigma$ through $q$ such that
	\begin{equation*}
		u(\sigma(t))=0,
		\qquad
		A(\sigma(t))\dot\sigma(t)=0
	\end{equation*}
	for all sufficiently small $|t|$.
\end{lemma}

\begin{proof}
	Since $q$ is a minimum point of $u$,
	\begin{equation*}
		\nabla^2u(q)\geq 0.
	\end{equation*}
	Thus $h_q\leq 0$. Together with $\det A(q)\leq 0$ and
	$A(q)\neq 0$, this shows that $A(q)$ has eigenvalues $0$ and
	a strictly negative number.

	After shrinking the neighborhood of $q$, we can choose smooth
	orthonormal eigenvector fields $T,E$ such that
	\begin{equation*}
		AT=\mu T,\qquad AE=-\eta E,
		\qquad
		\mu\geq 0,\quad \eta>0.
	\end{equation*}
	The nonnegativity of $\mu$ follows from $\det A\leq 0$.
	The Codazzi equation gives
	\begin{equation}
		\label{eq:three-eigenvector-codazzi}
		\nabla_TT
		=
		\frac{E\mu}{\mu+\eta}\,E.
	\end{equation}
	Indeed, comparing the $T$-components of
	$(\nabla_EA)T$ and $(\nabla_TA)E$ yields this identity.

	Let $\sigma$ be an integral curve of $T$ with $\sigma(0)=q$,
	and put
	$
	y(t)=u(\sigma(t)).
	$
	Using \eqref{eq:three-eigenvector-codazzi}, we obtain
	\begin{equation*}
		\begin{aligned}
		y''
		&=
		\nabla^2u(T,T)
		+
		du(\nabla_TT)\\
		&\leq
		-b\mu+\frac{E\mu}{\mu+\eta}\,Eu.
		\end{aligned}
	\end{equation*}

	On a sufficiently small relatively compact neighborhood,
	there are constants $b_0,\eta_0>0$ such that
	\begin{equation*}
		b\geq b_0,\qquad \eta\geq\eta_0.
	\end{equation*}
	The elementary local gradient estimate for nonnegative smooth
	functions, obtained from Taylor's formula, gives
	\begin{equation*}
		|\nabla u|^2\leq c\,u,
		\qquad
		|\nabla\mu|^2\leq c\,\mu.
	\end{equation*}
	Consequently, Young's inequality implies
	\begin{equation*}
		y''
		\leq
		-b_0\mu+c_1\sqrt{y\mu}
		\leq c_2y.
	\end{equation*}
	Moreover,
	\begin{equation*}
		y\geq 0,\qquad y(0)=y'(0)=0.
	\end{equation*}
	Integrating twice gives, for $t\geq 0$ sufficiently small,
	\begin{equation*}
		0\leq y(t)
		\leq
		c_2\int_0^t(t-s)y(s)\,ds.
	\end{equation*}
	Taking the maximum on a sufficiently short interval shows that
	$y\equiv 0$. The same argument applied to $y(-t)$ gives the
	conclusion for negative $t$.

	Along this segment, $du=0$, and hence
	\begin{equation*}
		0=y''
		=
		\nabla^2u(T,T)
		\leq -b\mu.
	\end{equation*}
	Thus $\mu=0$ along $\sigma$. Since $\mu\geq 0$ in the
	neighborhood, its gradient also vanishes along $\sigma$.
	Equation \eqref{eq:three-eigenvector-codazzi} now gives
	$
	\nabla_TT=0.
	$
	Therefore $\sigma$ is a unit-speed geodesic and
	$A\dot\sigma=0$ along it.
\end{proof}

\subsection{Rigidity from support-map fibers}
\label{sec:three-fiber-rigidity}

To choose a regular level meeting the contact set, we use the
following consequence of G\'alvez--Mira's global result on $\Sph^2$.

\begin{theorem}[Consequence of G\'alvez--Mira {\cite[Theorem~1.4 and Section~3]{GM}}]
	\label{thm:GM}
	Let $f\in C^2(\Sph^2)$ satisfy $\det A_f\leq0$ and
	$A_f\not\equiv0$. Then there are two closed geodesic semicircles
	$\Gamma_1,\Gamma_2\subset\Sph^2$ and distinct vectors
	$z_1,z_2\in\mathbb R^3$ such that the support map satisfies
	\begin{equation*}
		F_f(x)=\nabla f(x)+f(x)x=z_i
		\qquad\text{for every }x\in\Gamma_i,\quad i=1,2.
	\end{equation*}
\end{theorem}

The distinct values follow from the gradient-map statement in
\cite[Section~3]{GM}, applied to the homogeneous extension
$\widetilde f(rx)=rf(x)$, for which
$\nabla\widetilde f|_{\Sph^2}=F_f$.

\begin{proposition}
	\label{prop:three-contact-rigidity}
	Let $A$ be a smooth self-adjoint Codazzi tensor on $\Sph^2$
	with $\det A\leq0$, and put $h(X,Y)=\langle AX,Y\rangle$.
	Suppose that an open set $U\subset\Sph^2$ and functions
	$u,b\in C^\infty(U)$ satisfy
	\begin{equation*}
		u\geq0,\qquad b>0,\qquad \nabla^2u+b h\leq0
		\qquad\text{on }U,
	\end{equation*}
	and that $C=\{x\in U:u(x)=0\}$ is nonempty and compact.
	Then $A\equiv0$ on $\Sph^2$.
\end{proposition}

\begin{proof}
	At each point of $U$, the condition $\det A\leq0$ gives a
	nonnegative eigenvalue of $A$. Evaluating the tensor inequality
	on a corresponding unit eigenvector shows that
	\begin{equation*}
		\lambda_{\min}(\nabla^2u)\leq 0.
	\end{equation*}
	Lemma~\ref{lem:three-convex-barrier} therefore applies to $C$.
	In particular, every nonempty relatively clopen subset of $C$
	meets every closed geodesic semicircle.

	At every $q\in C$, the minimum property of $u$ also gives
	\begin{equation}
		\label{eq:three-contact-rank}
		h_q\leq 0,\qquad
		\det A(q)=0,\qquad
		\operatorname{rank}A(q)\leq 1.
	\end{equation}

	Suppose that $A\not\equiv0$, and let
	$F=\nabla f+fx$ be the support map from
	Lemma~\ref{lem:codazzi-potential-map}, so that
	\begin{equation}\label{eq:three-potential-map}
		dF=A.
	\end{equation}
	By Theorem~\ref{thm:GM}, we may choose two closed geodesic
	semicircles $\Gamma_1,\Gamma_2$ with
	\begin{equation*}
		F|_{\Gamma_i}=z_i,\qquad z_1\neq z_2.
	\end{equation*}

	Choose a linear functional $L:\mathbb R^3\to\mathbb R$ such that
	$
	a_1:=L(z_1)<L(z_2)=:a_2,
	$
	and put
	\begin{equation*}
		G=L\circ F:S^2\longrightarrow\mathbb R.
	\end{equation*}
	By Lemma~\ref{lem:three-convex-barrier}, each $\Gamma_i$ meets
	$C$, so $a_i\in G(C)$. In fact,
	\begin{equation*}
		[a_1,a_2]\subset G(C).
	\end{equation*}
	Otherwise, for some $s\in(a_1,a_2)\setminus G(C)$, the set
	\begin{equation*}
		K=C\cap\{G<s\}
	\end{equation*}
	would be a nonempty relatively clopen subset of $C$ disjoint
	from $\Gamma_2$, contradicting the same lemma.

	The contact inequality will now provide kernel propagation,
	even though $A$ need not have constant rank near $C$.
	By Sard's theorem, choose a regular value
	$
	s\in(a_1,a_2)
	$
	of $G$. Let $\Sigma$ be a connected component of $G^{-1}(s)$
	which meets $C$. Then $\Sigma$ is a compact connected embedded
	curve without boundary.

	For every $x\in\Sigma\cap C$, regularity gives
	\begin{equation*}
		0\neq dG_x=L\circ A(x).
	\end{equation*}
	Together with \eqref{eq:three-contact-rank}, this implies
	\begin{equation*}
		\operatorname{rank}A(x)=1,
		\qquad
		T_x\Sigma=\ker dG_x=\ker A(x).
	\end{equation*}
	Lemma~\ref{lem:three-rank-one-propagation} provides a geodesic
	segment through $x$ contained in $C$ and tangent to $\ker A$
	along the segment. By \eqref{eq:three-potential-map}, $F$,
	and hence $G$, is constant along this segment. It therefore
	lies in $\Sigma$ and forms a neighborhood of $x$ in $\Sigma$.

	Thus $\Sigma\cap C$ is relatively open in $\Sigma$. It is also
	closed and nonempty, so connectedness gives
	$
	\Sigma\subset C.
	$
	The same propagation argument at every point of $\Sigma$
	shows that $\Sigma$ is totally geodesic. Since it is compact,
	it is complete, and therefore it is a full great circle of $S^2$.

	Moreover, $T\Sigma=\ker A$ and $dF=A$, so $F$ is constant on
	$\Sigma$, say $F=z$. Thus $f(x)-\langle z,x\rangle=0$ on
	$\Sigma$. Lemma~\ref{lem:circle-support-rigidity} implies that
	this function is linear on $\Sph^2$. Its spherical Hessian is
	$A$, so $A\equiv0$, a contradiction.
\end{proof}

\begin{proof}[Proof of Theorem~\ref{thm:inradius-three}]
	Choose $p\in\operatorname{int}M$ such that
	$
	R=d_g(p,\partial M)=r_{\mathrm{in}}(M).
	$
	At any nearest boundary point $q$, Lemma~\ref{lem:three-distance}
	gives
	\begin{equation*}
		h_q\leq \cot R\,\bar{g}_q,
		\qquad R<\pi.
	\end{equation*}
	If $R>\pi/2$, then $h_q$ is negative definite, so
	$\det A(q)>0$. This contradicts \eqref{eq:three-gauss}.
	Therefore
	$
	r_{\mathrm{in}}(M)\leq\frac{\pi}{2}.
	$

	Assume now that $R=\pi/2$, and define the contact set
	\begin{equation*}
		C=\left\{q\in\partial M:
		d_g(p,q)=\frac{\pi}{2}\right\}.
	\end{equation*}
	This is a nonempty compact subset of $S^2$.

	By Lemma~\ref{lem:three-distance}, the function $r=d_g(p,\cdot)$
	is smooth near $C$. On a sufficiently small boundary
	neighborhood $U$ of $C$, set
	\begin{equation*}
		u=r|_{\partial M}-\frac{\pi}{2},
		\qquad
		b=\langle\nabla^M r,\nu\rangle.
	\end{equation*}
	We can arrange that
	\begin{equation*}
		u\geq 0,\qquad
		C=\{u=0\},\qquad
		b>0,\qquad
		\frac{\pi}{2}\leq r<\pi
	\end{equation*}
	on $U$.

	The nearby minimizing geodesics lie in $M$, as in
	Lemma~\ref{lem:three-distance}, so the boundary Hessian formula
	and Hessian comparison under $\sec_g\geq1$ give
	\begin{equation*}
		\begin{aligned}
		\nabla^2u+b h
		&=
		\nabla_M^2r|_{T\partial M}\\
		&\leq
		\cot r\,(\bar{g}-du\otimes du)
		\leq 0.
		\end{aligned}
	\end{equation*}
	Here $|du|_{\bar{g}}\leq 1$, so $\bar{g}-du\otimes du$ is
	positive semidefinite, whereas $\cot r\leq 0$.
	Proposition~\ref{prop:three-contact-rigidity} now gives $A\equiv0$,
	so $\partial M$ is totally geodesic.
	Since $\sec_g\geq1$ implies $\Ric_g\geq2g$, the hemisphere
	rigidity theorem of Hang--Wang~\cite{HW} gives
	$
		(M^3,g)\cong S^3_+.
	$
	Conversely, the standard upper hemisphere has inradius
	$\pi/2$. This completes the proof.
\end{proof}

\medskip
\noindent\textbf{Declaration of AI use.}
We used ChatGPT for proof exploration and editing.
We take responsibility for verification and exposition.

\begingroup
\renewcommand*{\bibfont}{\small}
\setlength{\emergencystretch}{2em}
\printbibliography
\endgroup
\end{document}